\documentclass[a4paper, 12pt]{article}
\usepackage[a4paper]{geometry}
\usepackage{amsmath}
\usepackage{amssymb}
\usepackage{latexsym}
\usepackage{enumerate}
\usepackage{graphics}
\usepackage[dvips]{graphicx}
\usepackage{epsfig}
\usepackage{lineno}
\usepackage{enumerate}
\usepackage[english]{babel}
\usepackage{amsmath,amssymb,amsthm}
\usepackage{subfigure}
\newcommand{\be}{\begin{equation}}
\newcommand{\ee}{\end{equation}}
\newcommand{\bea}{\begin{eqnarray}}
\newcommand{\eea}{\end{eqnarray}}

\newcommand{\nod}{\noindent}

\newcommand{\ba}{\begin{array}}
\newcommand{\ea}{\end{array}}
\newcommand{\bc}{\begin{center}}
\newcommand{\ec}{\end{center}}

\allowdisplaybreaks

\begin{document}
%\linenumbers*[1]
%\title{\bf Exact and approximate epidemic models on networks: a new, improved closure relation\\}
\title{\bf A class of pairwise models for epidemic dynamics on weighted networks\\}

\author{Prapanporn Rattana$^{1}$, Konstantin B. Blyuss $^{1}$, Ken T.D. Eames$^{2}$ \& Istvan Z. Kiss$^{1,\ast}$}

\maketitle

\begin{center}

$^1$School of Mathematical and Physical Sciences, Department of
Mathematics, University of Sussex, Falmer,
Brighton BN1 9QH, UK\\
$^2$The Centre for the Mathematical Modelling of Infectious Diseases, London School of Hygiene and Tropical Medicine, Keppel Street, London WC1E 7HT\\

\end{center}

\vspace{12cm}
\begin{flushleft}
$\ast$corresponding author\\
email: i.z.kiss@sussex.ac.uk\\
\end{flushleft}

\newpage
\begin{abstract}
In this paper, we study the $SIS$ (susceptible-infected-susceptible) and $SIR$ (susceptible-infected-removed) epidemic models
on undirected, weighted networks by deriving pairwise-type approximate models coupled with individual-based network simulation.
Two different types of theoretical/synthetic weighted network models are considered. Both models start from non-weighted networks
with fixed topology followed by the allocation of link weights in either (i) random or (ii) fixed/deterministic way. The pairwise models
are formulated for a general discrete distribution of weights, and these models are then used in conjunction with network simulation
to evaluate the impact of different weight distributions on epidemic threshold and dynamics in general. For the $SIR$ dynamics, the
basic reproductive ratio $R_0$ is computed, and we show that (i) for both network models $R_{0}$ is maximised if all weights are equal,
and (ii) when the two models are equally matched, the networks with a random weight distribution give rise to a higher $R_0$ value. The
models are also used to explore the agreement between the pairwise and simulation models for different parameter combinations.
\end{abstract}

\nod {\bf Keywords:} weighted network, pairwise model, moment closure

\newpage

%%%%%%%%%%%%%%%%%%%%%%%%%%%%%%%%%%%%%%%%%%%%%%%%%%%%%%%%%%%%%%%%%%%%%%%%%%%%%%%%%%%%%%%%%%
\section{Introduction}
%%%%%%%%%%%%%%%%%%%%%%%%%%%%%%%%%%%%%%%%%%%%%%%%%%%%%%%%%%%%%%%%%%%%%%%%%%%%%%%%%%%%%%%%%%

Conventional models of epidemic spread consider a host population of identical individuals, each interacting in the same way with each of the others (see \cite{AndersonMay,DiHe,KeRo} and references therein). At the same time, in order to develop more realistic mathematical models for the spread of infectious diseases, it is important to obtain the best possible representation of the corresponding transmission mechanism. To achieve this, more recent models have included some of the many complexities that have been observed in mixing patterns. One such approach consists in splitting the population into a set of different subgroups, each with different social behaviours. Even more detail is included within network approaches, which allow to include differences between individuals, not just between sub-populations. In such models, each individual is represented as a node, and interactions that could permit the transmission of infection appear as edges linking nodes. The last decade has seen a substantial increase in the research of how infectious diseases are spread over large networks of connected nodes, where networks themselves can represent either small social contact networks or larger scale travel networks, including global aviation networks \cite{DanonNetwReview,DM,GS,NetwReviewKeelingEames,MoNe,MPV,NewEpi,PS1,PS2}. Importantly, the characteristics of the network, such as the average degree and the node degree distribution have a profound effect on the dynamics of the infectious disease spread, and hence significant efforts are made to capture properties of realistic contact networks.

Network models provide an intuitively appealing way to consider the social structure of populations, but they bring with them the challenge of collecting sufficient data to parameterise them fully. Ideally, an entire population would be sampled and all relevant interactions measured in order to reconstruct the contact network, but this is, unsurprisingly, seldom possible, and it is generally necessary to make suitable approximations. With the advance of computer-based tracking, statistical properties of many realistic networks (e.g. mobile phone calls) have been studied in great detail, and these are often used as proxy for the analysis of epidemic dynamics \cite{DM,HBG}.

One of the simplifying assumptions often put into models on networks is that all links are equally likely to transmit infection \cite{Boc,Eubank,NetwReviewKeelingEames,Riley07}. However, a more detailed consideration leads to an observation that this is often not the case, as some links are likely to be far more capable of transmitting infection than others due to closer contacts (e.g. within households \cite{Beu}) or long-duration interactions \cite{Coh97,Ed97,Read08,Riley07,RFer06}.
%For example, living in the same house as an infected individual is likely to be riskier than travelling in the same train with them.
To account for this heterogeneity in properties of social interactions, network models can be adapted,
%by assigning weights to the network links
thus resulting in {\it weighted contact networks}, where connections between different nodes have different weights. These weights may be associated with the duration,
proximity, or social setting of the interaction, and the key point is that they are expected to be correlated with the risk of disease transmission. The precise relationship between the properties of an interaction and its riskiness is hugely complex; here, we will consider a``weight" that is exactly proportional to the transmission rate
along a link. From a general perspective, one can consider two types of weighted networks: symmetric networks, where the strength only varies between different edges but is the same for a given pair of connected nodes, or go one step further and allow for the strength to be different even for the same pair of connected nodes $A$ and $B$, depending on whether $A$ infects $B$, or $B$ infects $A$. The first type of networks can be efficiently conceived as undirected graphs represented by symmetric connection matrices, while for the second one it is necessary to consider the contact network as a directed graph, resulting in a connection matrix, which in general is non-symmetric.
%In the case when an epidemic is considered to be spreading on a network of particular type, it is often assumed that the strength of connection is the same for all nodes. However, a more realistic representation of contact networks can be yielded by using {\it weighted contact networks}, where connections between different nodes have different weights.
Although consideration of weighted networks may seem as an additional complication for the analysis of epidemic dynamics, in fact it provides a much more realistic representation of actual contact networks.

Early models of weighted networks arose in the context of identifying genetic and metabolic networks \cite{Al,NS}. Later on, they found applications in a number of diverse research areas, such as
international trade networks \cite{BMSKM,FRS}, scientific collaboration networks \cite{BBPV,Fan,LiWu}, aviation networks of individual countries \cite{Bag,LiCai} or global aviation networks 
\cite{r1,Col2,HBG} and networks of mobile phone communications \cite{OS}. Yook {\it et al.} \cite{YJBT} started a systematic study of evolving networks with nonbinary connectivities by considering growing weighted and exponential networks. This approach was later developed into a generic formalism for the analysis of evolving weighted networks \cite{r1,r7,r8,LC,LiWu}. Newman \cite{New} has
discussed general properties of weighted networks. Substantial amount of work has been done on the analysis of scale-free networks with different types of weight distribution \cite{WZ,ZTZH}. Wang {\it
et al.} \cite{WW} have proposed a model with dynamical adjustment of weights in technological networks.

In epidemiological context, Yan {\it et al.} \cite{Yan} analysed a
model on weighted scale-free networks and found that heterogeneity
in weight distribution leads to a slowdown in the spread of
epidemics. Furthermore, they have shown that for a  given network
topology and mean infectivity, epidemics spread fastest in
unweighted networks. Chu {\it et al.} \cite{r17} have investigated
the dynamics of an $SI$ model on weighted scale-free networks with
community structure and showed that the contribution coming from
weights of edges connecting different communities exceeds the
impact of weights on edges inside communities. Karsai {\it et al.}
\cite{Ka} have studied mean-field and finite-size scaling in
weighted networks with small weights on edges connecting
high-degree nodes. Yang {\it et al.} \cite{Yang} have shown that
disease prevalence can be maximized when the edge weights are
chosen to be inversely proportional to the degrees of receiving
nodes but, in this case, the transmissibility was not directly
proportional to the weights and weights were also asymmetric. Yang
\& Zhou \cite{r16} have considered $SIS$ epidemics on homogeneous
networks with uniform or power-law edge weight distribution and
shown how to derive a certain type of mean-field description for
such models. Britton {\it et al.} \cite{BDL} have derived an
expression for the basic reproductive ratio in weighted networks
with generic distributions of node degree and weights, and Deijfen
\cite {Deij} has performed a similar analysis to study
vaccination in such networks. In terms of practical
epidemiological applications, weighted networks have already been
effectively used to study control of global pandemics
\cite{Col,Coop,r12} and the spread of animal disease due to cattle
movement between farms \cite{GMB,Webb}. Eames {\it et al.}
\cite{r12} have considered an $SIR$ model on an undirected
weighted network, where rather than using some theoretical
formalism to generate an idealized network, the authors have
relied on social mixing data obtained from questionnaires
completed by members of a peer group \cite{Read08} to construct a
realistic weighted network. Having analysed the dynamics of
epidemic spread in a such a network, they showed how  information
about node-specific infection risk can be used to develop targeted
preventative vaccination strategies. As an alternative approach to
modelling heterogeneity in network interactions, Eames \cite
{Eames2008} has considered a network with two types of contacts:
random and regular, where the first type refers to arbitrary
connections with any node in the network, while the second type
designates multiple repeated contacts with the same nodes. He
showed that in a highly clustered population, random contacts
allow infection to reach otherwise inaccessible parts of the
network, while in the case of all contacts being regular,
clustering leads to a significant reduction in the spread of
infection.

In this paper, we consider the dynamics of an infectious disease spreading on weighted networks
with different weight distributions. Since we are primarily concerned with the effects
of weight distribution on the disease dynamics, the connection matrix will be assumed to be symmetric,
representing the situation when the weights can only be different for different network edges, but for a
given edge the weight is the same irrespective of the direction of infection. From epidemiological perspective,
we consider both the case when the disease confers permanent immunity (represented by an $SIR$ model),
and the case when the immunity is short-lived, and upon recovery the individuals return to the class of
susceptibles ($SIS$ model). For both of these cases we derive the corresponding ODE-based pairwise
models and their closure approximations. Numerical simulation of both the network dynamics and the
pairwise approximations are performed.

The outline of this paper is as follows. In the next section, the construction of specific weighted networks to be used for
analysis of epidemic dynamics is discussed. This is complemented by the derivation of corresponding pairwise models and their
closure approximations. Section 3 contains the derivation of the basic reproductive ratio $R_0$ for the $SIR$ model and for different
weight distributions as well as numerical simulation of both network
models and their pairwise ODE counterparts. The paper concludes in
Section 4 with discussion of results and possible further extensions of this work.
%%%%%%%%%%%%%%%%%%%%%%%%%%%%%%%%%%%%%%%%%%%%%%%%%%%%%%%%%%%%%%%%%%%%%%%%%%%%%%%%%%%%%%%%%%
\section{Material and methods}
%%%%%%%%%%%%%%%%%%%%%%%%%%%%%%%%%%%%%%%%%%%%%%%%%%%%%%%%%%%%%%%%%%%%%%%%%%%%%%%%%%%%%%%%%%

%%%%%%%%%%%%%%%%%
\subsection{Network construction}
%%%%%%%%%%%%%%%%%

There are two conceptually different approaches to constructing
weighted networks for modelling infectious disease spread. In the first approach, there is a seed or a primitive
motif, and the network is then grown or evolved from this initial
seed according to some specific rules. In this method, the
topology of the network is co-evolving with the distribution of
weights on the edges \cite {r7,r8,r9,LC,Yang}. Another approach is
to consider a weighted network as a superposition of an unweighted
network with a distribution of weights across edges which could
independent of the original network or it may be correlated with
node metrics, such as their degree, \cite{BDL,Deij,r11,r3}. In this
paper we use the second approach in order to investigate the
particular role played by the distribution of weights across
edges, rather than network topology, in the dynamics of epidemic
spread. Besides computational efficiency, this will allow us to
make some analytical headway in deriving and analysing
low-dimensional pairwise models which are likely to perform better
when weights are attached according to the scenarios described
above.

%In the simplest case, a weighted graph or network can be considered as the superposition of an unweighted network with a prescribed
%degree distribution and a given weight distribution across edges. Weights can be distributed at random across all edges as given by the weight
%distribution or to achieve some form of correlation between node degree and weights or between the weights of the links belonging to the same node, i.e. correlated link weights around nodes.
%\textbf{Need some references for different types of weighted networks!\\}

Here we consider two different methods of assigning weights to network links: a network in which weights are assigned to links at random, and
a network in which each node has the same distribution of weighted links connected to it. In reality, there is likely to be a great deal more structure
to interaction weights, but in the absence of precise data and also for the purposes of developing models that allow one to explore a number of different
assumptions, we make these simplifying approximations.

%Due to this added level of complexity many scenarios can be considered, but here, we will consider some of the simpler cases.
%This choice is motivated by our aim to \\

%%%%%%%%%
\subsubsection{Random weight distribution}
%%%%%%%%%
First we consider a simple model of an undirected weighted network with $N$ nodes where the weights of the links can take values
$w_i$ with probability $p_i$, where $i=1,2, \dots, M$. The underlying degree distribution of the corresponding unweighted
network can be chosen to be of the more basic forms, e.g. homogeneous random or Erd\H{o}s-R\'enyi-type random networks.

The generation of such networks is straightforward, and weights can be assigned during link creation in the unweighted network.
For example, upon using the configuration model for generating unweighted networks, each new link will have a weight assigned to it based
on the chosen weight distribution. This means that in a homogeneous random network with each node having $k$ links, the distribution of link
weights of different type will be multinomial, and it is given by
\begin{equation}
P(n_{w_1}, n_{w_2}, \dots, n_{w_M})=\frac{k!}{n_{w_1}! n_{w_2}! \dots n_{w_M}!}p_1^{n_1}p_2^{n_2}\dots p_M^{n_M},
\end{equation}
where, $n_{w_1}+n_{w_2}+ \cdots+ n_{w_M}=k$ and $P(n_{w_1}, n_{w_2}, \dots, n_{w_M})$ stands for the probability of a node having $n_{w_1}$, $n_{w_2}$, $\dots$, $n_{w_M}$ links with weights $w_1$, $w_2$, $\dots$, $w_M$, respectively.
%The calculation above becomes more complicated if Erd\H{o}s-R\'enyi random networks with heterogenous degree distributions are considered.
While the above expression is applicable in the most general set-up, it is worth considering the case of weights of only two types, where
the distribution of link weights for a homogenous random network becomes binomial
\begin{equation}
P(n_{w_1}, n_{w_2}=k-n_{w_1})=
{k \choose n_{w_1}}p_1^{n_1}(1-p_1)^{k-n_1},
\end{equation}
where, $p_1+p_2=1$ and $n_{w_1}+n_{w_2}=k$.  The average link weight in the model above can be easily found as
\[
w_{av}^{random}=\sum_{i=1}^{M}p_iw_i,
\]
which for the case of weights of two types $w_1$ and $w_2$ reduces to
\[
w_{av}^{(2r)}=p_1w_1+p_2 w_2=p_1w_1+(1-p_1)w_2.
\]

%Similarly for Erd\"os-R\'enyi random networks, the distribution above becomes
%\begin{equation}
%P(n_{w_1}, n_{w_2}=k-n_{w_1})=\sum_{i=n_{w_1}}^{N-1}
%{i \choose n_{w_1}}P(w_1)^{n_1}(1-P(w_1))^{i-n_1}P_{ER}(i),
%\end{equation}
%\textbf{I don't think the above equation is right. No need fo the sum, and the i should be a k)}
%where $P_{ER}={N-1 \choose i}p^i(1-p)^{N-1-i}$ is the probability
%of a node in an Erd\"os-R\'enyi random network, with an average
%degree $\langle k \rangle=p(N-1)$, having degree $i$.

%%%%%%%%%
\subsubsection{Fixed deterministic weight distribution}
%%%%%%%%%
%\textbf{Ken, would you please add your bit here, i.e. networks with fixed degree $k$, but with say $k_1$ and $k_2$ links with weights $w_1$ and $w_2$, respectively.
%Some motivation would also be good. I assume that you have in mind that almost all nodes have a number of strong connections with say family plus a few weak links or occasional contacts.}

As a second example we consider a network, in which each node has $k_i$ links with weight $w_i$ ($i=1, 2, \dots, M$), where $k_1+k_2+ \cdots+ k_M=k$. The different weights here could be interpreted as  being associated with different types of social interaction: e.g. home, workplace, and leisure contacts, or physical and non-physical interactions. In this model all individuals are identical in terms of their connections, not only having the same number of links (as in the model above) but also having the same set of weights. The average weight in such a model is given by
\[
w_{av}^{fixed}=\sum_{i=1}^{M}p_iw_i,\hspace{0.5cm}p_i=\frac{k_i}{k},
\]
where $p_i$ is the fraction of links of type $i$ for each node. In the case of links of two types with weights $w_1$ and $w_2$, the average weight becomes
\[
w_{av}^{(2f)}=p_1w_1+p_2 w_2=\frac{k_1}{k}w_1+\frac{k_2}{k}w_2=\frac{k_1}{k}w_1+\frac{k-k_1}{k}w_2.
\]

%%%%%%%%%
\subsubsection{Epidemic models}
%%%%%%%%%
In this study, the simple $SIS$ and $SIR$ epidemic models are considered.
The epidemic dynamics is specified in terms of infection and recovery events. The rate of transmission
across an unweighted edge between an infected and susceptible individual is denoted by $\tau$. This will then be
adjusted by the weight of the link which is assumed to be directly proportional to the strength of the transmission along that link.
Infected individuals recover independently of each other at rate $\gamma$. The simulation is implemented using
Gillespie algorithm \cite{Gil} with inter-event times distributed exponentially with a rate given by the total rate of change in the network, with the single
event to be implemented at each step being chosen at random and proportionally to its rate. All simulations start with most nodes being susceptible
and with a few infected nodes chosen at random.

%%%%%%%%%%%%%%%%%
\subsection{Pairwise equations}
%%%%%%%%%%%%%%%%%
In this section we extend the classic pairwise model for unweighted networks \cite{Keeling1999,r5} to the case of weighted graphs with $M$ different link-weight types.
Pairwise models successfully interpolate between classic compartmental ODE models and full individual-based network simulation with the added advantage of high
transparency and a good degree of analytical tractability. These qualities makes them an ideal tool for studying dynamical processes on
networks \cite{Eames2008,Vasilis,ThomasUnifying,Keeling1999}, and they can be used on their own and/or in parallel with simulation.
The original versions of the pairwise models have been successfully extended to networks with heterogenous degree distribution \cite{HetPairwise}, asymmetric networks \cite{KieranAsymmetric} and
situations where transmission happens across different/combined routes \cite{Eames2008,Vasilis} as well as when taking into consideration network motifs of higher order than pairs and triangles \cite{HouseMotif}. The extension that we propose is based on the previously established precise counting procedure at the
level of individuals, pairs and triples, as well as on a careful and systematic account of all possible transitions needed to derive the full set of evolution equations for singles and pairs. These obviously
involve the precise dependency of lower order moments on higher order ones, e.g. the rate of change of the expected number of susceptible nodes is proportional to the
expected number of links between a susceptible and infected node. We extend the previously well-established notation \cite{Keeling1999} to account for the added level
of complexity due to different link weights. In line with this, the number of singles remains unchanged, with $[A]$ denoting the number of nodes across
the whole network in state $A$. Pairs of type $A-B$, $[AB]$, are now broken down depending on link weights, i.e. $[AB]_{i}$ represents the number of links of
type $A-B$ with the link having weight $w_i$, where as before $i=1, 2, \dots, M$ and $A, B \in \{S, I, R\}$ if an $SIR$ dynamics is used. As before, links are
doubly counted (e.g. in both directions) and thus the following relations hold: $[AB]_m=[BA]_m$ and $[AA]_m$ is equal to twice the number of uniquely counted
links of weight $w_m$ with nodes at both ends in state $A$. From this extension it follows that $\sum_{i=1}^{M} [AB]_i=[AB]$. The same convention holds at the level of triples where
$[ABC]_{mn}$ stands for the expected number of triples where a node in state $B$ connects a node in state $A$ and $C$ via links of weight $w_{m}$ and $w_{n}$, respectively. The weight of the link
impacts on the rate of transmission across that link, and this is achieved by using a link-specific transmission rate equal to $\tau w_{i}$, where $i=1, 2, \dots, M$.
In line with the above, we construct two pairwise models,  one for $SIS$ and one for $SIR$ dynamics.

%\subsubsection{Evolution equations for the $SIS$ model}
The pairwise model for the $SIS$ dynamics can be written in the form:
\begin{equation}\label{SIS}
\begin{array}{l}
[\dot{S}]=\gamma{[I]}-\tau{\sum_{n=1}^M {w_{n}}{[SI]_{n}}},\\\\
\ [\dot{I}]=\tau{\sum_{n=1}^M {w_{n}}{[SI]_{n}}}-\gamma{[I]},\\\\
\ [\dot{SI}]_{m}= \gamma{([II]_{m}-[SI]_{m})}+\tau{\sum_{n=1}^M {w_{n}}{([SSI]_{mn}-[ISI]_{nm})}}-\tau{w_{m}}{[SI]_{m}},\\\\
\ [\dot{II}]_{m}= -2{\gamma{[II]_{m}}}+2{\tau}{\sum_{n=1}^M w_{n}{[ISI]_{nm}}}+2{\tau}{w_{m}}{[SI]_{m}},\\\\
\ [\dot{SS}]_{m}= 2{\gamma{[SI]_{m}}}-2{\tau{\sum_{n=1}^M w_{n}{[SSI]_{mn}}}},
\end{array}
\end{equation}
where $m=1,2,3,...,M$ and $[AB]_{m}$ denotes the expected number of links with weight $w_{m}$ connecting two nodes of type $A$ and $B$, respectively ($A,B\in \{S, I\}$).

%\subsubsection{Evolution equations for the $SIR$ model}
In the case when upon infection individuals recover at rate $\gamma$ and once recovered they maintain a life-long immunity, we have the following system of equations
describing the dynamics of a pairwise $SIR$ model:

\begin{equation}\label{SIR}
\begin{array}{l}
   \dot{[S]}= -\tau{\sum_{n=1}^M {w_{n}}{[SI]_{n}}},\\\\
   \dot{[I]}= \tau{\sum_{n=1}^M {w_{n}}{[SI]_{n}}}-\gamma{[I]},\\\\
   \dot{[R]}= \gamma{[I]},\\\\
   \ [\dot{S}{S}]_{m}= -2{\tau{\sum_{n=1}^M {w_{n}}{[SSI]_{mn}}}},\\\\
   \ [\dot{S}{I}]_{m}= \tau{\sum_{n=1}^M {w_{n}}{([SSI]_{mn}-[ISI]_{nm})}}-\tau{w_{m}}{[SI]_{m}}-\gamma{[SI]_{m}},\\\\
   \ [\dot{S}{R}]_{m}= -{\tau{\sum_{n=1}^M {w_{n}}{[ISR]_{nm}}}}+\gamma{[SI]_{m}},\\\\
   \ [\dot{I}{I}]_{m}= 2{\tau}{\sum_{n=1}^M
   w_{n}{[ISI]_{nm}}}+2{\tau}{w_{m}{[SI]_{m}}}-2{\gamma}{[II]_{m}},\\\\
   \ [\dot{I}{R}]_{m}= {\tau{\sum_{n=1}^M {w_{n}}{[ISR]_{nm}}}}+\gamma{([II]_{m}-[IR]_{m})},\\\\
   \ [\dot{R}{R}]_{m}= \gamma{[IR]_{m}},
   \end{array}
\end{equation}
where again $m=1,2,3,...,M$ with the same notation as above.

The above systems of equations (\ref{SIS}) and (\ref{SIR}) are not closed, as  equations for the pairs require knowledge of triples, and thus, equations
for triples are needed.  This dependency on higher-order moments can be curtailed by closing the equations via approximating triples in terms of singles and pairs \cite{Keeling1999}. For both systems, the agreement with simulation will heavily depend on the precise distribution of weights across the links, the network topology, and the type of closures that will be used to capture essential features of network structure and the weight distribution. As a check and a reference back to previous model, in Appendix A shows how systems (\ref{SIS}) and (\ref{SIR}) reduce to the standard unweighted pairwise $SIS$ and $SIR$ models \cite{Keeling1999} when all weights are equal to each other, $w_1=w_2=\cdots=w_M=W$.

%%%%%%%%%%%%%%%%%%%%%%%
\subsubsection{Closure relations}
%%%%%%%%%%%%%%%%%%%%%%%
The most natural extension of the classic closure is given by
\begin{equation}\label{Clos_cl}
[ABC]_{mn}=\frac{k-1}{k}\frac{[AB]_m[BC]_n}{[B]},
\end{equation}
where $k$ is the number of links per node for a homogeneous natwork or the average nodal degree for networks with other than homogenous degree distributions.
Even for the simplest case of homogenous random networks with two weights, a different version of the closure can be derived. The starting point for this different
closure is the observation that the average number of links of weight $w_1$ across the whole network is $k_1=p_1k \le k$, and similarly,  the average number of links
of weight $w_2$ is $k_2=(1-p_1)k \le k$. This motivates the following closures
\begin{equation}\label{Clos}
\begin{array}{l}
\displaystyle{\ [ABC]_{11}=[AB]_1(k_1-1)\frac{[BC]_1}{k_1[B]}=\frac{k_1-1}{k_1}\frac{[AB]_1[BC]_1}{[B]},}\\ \\
\displaystyle{\ [ABC]_{12}=[AB]_1 k_2 \frac{[BC]_2}{k_2[B]}=\frac{[AB]_1[BC]_2}{[B]},}\\\\
\displaystyle{\ [ABC]_{21}=[AB]_2 k_1 \frac{[BC]_1}{k_1[B]}=\frac{[AB]_2[BC]_1}{[B]},}\\\\
\displaystyle{\ [ABC]_{22}=[AB]_2(k_2-1)\frac{[BC]_2}{k_2[B]}=\frac{k_2-1}{k_2}\frac{[AB]_2[BC]_2}{[B]}.}
 \end{array}
 \end{equation}
The specific choice of closure will depend on the structure of the network and, especially, how the weights are distributed. For example, for the case of
the homogeneous random networks with links  allocate randomly, both closures offer a viable alternative. For the case of a network where each node has a fixed
pre-allocated number of links with different weights, e.g. $k_1$ and $k_2$ links with  weights $w_1$ and $w_2$, respectively, the second closure (\ref{Clos})
offers the more natural/intuitive avenue towards closing the system and obtaining good agreement with network simulation.

%%%%%%%%%%%%%%%%%%%%%%%%%%%%%%%%%%%%%%%%%%%%%%%%%%%%%%%%%%%%%%%%%%%%%%%%%%%%%%%%%%%%%%%%%%
\section{Results}
%%%%%%%%%%%%%%%%%%%%%%%%%%%%%%%%%%%%%%%%%%%%%%%%%%%%%%%%%%%%%%%%%%%%%%%%%%%%%%%%%%%%%%%%%%
%The bulk of our results are derived for the case of two distinct weights and include comparisons to the equal weights (or unweighted) case.

In this section we present analytical and numerical results for weighted networks and pairwise representations of $SIS$ and $SIR$
models in the case of two different link-weight types (i.e. $w_1$ and $w_2$).

%%%%%%%
\subsection{Threshold dynamics for the $SIR$ model - the network perspective}
%%%%%%%
%\textbf{Need to define $R_0$ and other unexplained notation below
%I don't find the next section very easy to read. Expressions need to be explained, I think. It would perhaps work better if both models were dealt with using the NGM rather than arguing heuristically.
%Entries in NGM need to be defined clearly. The hand-waving approximation to the model 2 $R_0$ might be worth including for a population with more than 2 weights, as an easily understood analytic
%approx.}
The basic reproductive ratio, $R_0$ (the average number of secondary cases produced by a typical index case in an otherwise susceptible population), is one of the most fundamental quantities in  epidemiology (\cite{AndersonMay,Diekmann1990R0}). Besides informing us on whether a particular disease will spread in a population,
as well as quantifying the severity of an epidemic outbreak, it can be also used to calculate a number of other important quantities that have good intuitive interpretation. In what follows, we will compute $R_0$ and $R_0$-like quantities and will discuss their relation to each other,  and also issues around these being
model-dependent.  First, we compute $R_0$ from an individual-based or network perspective by employing the next generation matrix approach as used in the context of models with multiple transmission routes such as household models \cite{BallOnKiss}.\\

\noindent \textit{Random weight distribution:} First we derive an expression for $R_0$ when the underlying network is homogeneous, and the weights of the links are assigned at random according to a prescribed weight distribution. In the spirit of the proposed approach, the next generation matrix can be easily computed to yield
\[ NGM= (a_{ij})_{i,j=1,2}=\left| \begin{array}{cc}
(k-1)p_1r_1& (k-1)p_1r_1\\
(k-1)p_2r_2 & (k-1)p_2r_2 \end{array} \right|,\]
where
\[
r_1=\frac{\tau w_1}{\tau w_1+\gamma},\hspace{0.5cm}r_2=\frac{\tau w_2}{\tau w_2+\gamma}
\]
represent the probability of transmission from an infected to a susceptible across a link of weight $w_1$ and $w_2$, respectively.
Here, the entry $a_{ij}$ stands for the average number of infections produced via links of type $i$ (i.e. with weight $w_i$) by a typical infectious node who itself has been infected across a link of type $j$
(i.e. with weight $w_j$). Using the fact that $p_2=1-p_1$, the basic reproductive ratio can be found from the leading eigenvalue of the $NGM$ matrix as follows
\begin{equation}\label{R01}
R_{0}^1=(k-1)(p_1r_1+(1-p_1)r_2).
\end{equation}
In fact, the expression for $R_{0}$ can be simply generalised to more than two weights to give $R_{0}=(k-1)\sum_{i=1}^{M}p_ir_i$, where $w_m$ has frequency given by $p_m$ with the constraint that $\sum_
{i=1}^{M}p_i=1$. It is straightforward to show that upon assuming uniform weight distribution $w_i=W$ for $i=1,2, \dots, M$, the basic reproduction number on a homogeneous graph reduces to $R_0=(k-1)r$ as expected, and where, $r=\tau W/(\tau W+\gamma)$.\\

%In deriving $R_{0}$ we will show that due to the weights of the links being assigned at random, it is not necessary to consider whether infection from an initial infected node chosen uniformly at
%random spreads across a link of weight $w_1$ or $w_2$. This observation simplifies the expression for the basic reproduction number which can be written as:
%\begin{eqnarray}
%R_{0}&=&P^{w_1}\sum_{i=0}^{k}\frac{iP^{w_1}(i)}{\sum_{j=0}^{k}jP^{w_1}(j)}((i-1)r_1+(k-i) r_2)\\ \notag
%&+&P^{w_2}\sum_{i=0}^{k}\frac{iP^{w_2}(i)}{\sum_{j=0}^{k}jP^{w_2}(j)}((k-i)r_1+(i-1) r_2)\\ \notag
%&=&(k-1)(p_1r_1+(1-p_1)r_2),
%\end{eqnarray}
%where $r_1=\tau w_1/(\tau w_1+\gamma)$ and $r_2=\tau w_2/(\tau w_2+\gamma)$ represent the probability of transmission from an infected to a susceptible across a link of weight $w_1$ and
%$w_2$, respectively. The probability that the first infection from a randomly chosen index case has spread along a link of weight $w_1$ and $w_2$ is denoted by
%$$P^{w_1}=\sum_{i=0}^{k}{k \choose i}P(w_1)^{i}(1-P(w_1))^{k-i}\frac{ i w_1}{i w_1+(k-i) w_2},$$
%and
%$$P^{w_2}=\sum_{i=0}^{k}{k \choose i}P(w_1)^{i}(1-P(w_1))^{k-i}\frac{(k-i) w_2}{i w_1+(k-i) w_2}=1-P^{w_1}.$$
%Furthermore, $P^{w_2}(i)=P^{w_1}(k-i)$ and $\frac{iP^{w_1}(i)}{\sum_{j=0}^{k}jP^{w_1}(j)}$ is the excess degree of links of weight $w_1$ when transmission has reached a node across a link of
%weight $w_1$.

\noindent \textit{Deterministic weight distribution:} The case when the number of links with given weights for each node is fixed can be captured with the same approach, and the next generation matrix can be constructed as follows
\[ NGM= \left| \begin{array}{cc}
(k_1-1)r_1& k_1r_1\\
k_2r_2 & (k_2-1)r_2 \end{array} \right|.\]
As before, the leading eigenvalue of the $NGM$ matrix yields the basic reproductive ratio,
\begin{equation}\label{R02}
R_{0}^2=\frac{(k_1-1)r_1+(k_2-1)r_2+\sqrt{[(k_1-1)r_1-(k_2-1)r_2]^2+4k_1k_2r_1r_2}}{2}.
\end{equation}
%It is worth mentioning that the transmission potentials, i.e. the number of secondary infections produced by a randomly chosen single index case, for the two scenarios are
%$\rho_{0}^{1}=\sum_{i=0}^{k}(i r_1+(k-i) r_2)P^{w_1}(i)$ and $\rho_{0}^{2}=k_1r_1+k_2 r_2$, respectively.\\
%\textbf{We can use these expressions to say something about the role of weight distribution? Do we want to expand on this?\\}

\noindent Using these two equivalent expressions for the basic reproductive ratio, it is possible to prove the following result.\\

%\begin{theo} 
\noindent {\bf Theorem 1.} {\it Given the setup for the fixed weight distribution and using $p_1=k_1/k$, $p_2=k_2/k$ and $k_1+k_2=k$, if $1\leq k_1\leq k-1$ (which implies that $1\leq k_2\leq k-1$), then $R_{0}^{2} \le R_{0}^{1}$.}\\

\noindent The proof of this result is sketched out in Appendix B. This Theorem effectively states that provided each node has at least one link of type 1 and one link of type 2, then independently of disease parameters, it follows that the basic reproductive ratio as computed from (\ref{R01}) always exceeds or is equal to an equivalent $R_0$ computed from (\ref{R02}).

It is worth noting that both $R_0$ values reduce to
\begin{equation}
R_{0}^1=R_{0}^2=R_{0}=(k-1)r=\frac{(k-1)\tau W}{\tau W+\gamma},
\end{equation}
if one assumes that weights are equal, i.e. $w_1=w_2=W$. As one would expect, the first good indicator of the impact of weights on the epidemic dynamics will be the
average weight. Hence, it is worth considering the problem of maximising the values $R_0$ under assumption of a fixed average weight:
\begin{equation}
p_1w_1+p_2w_2=W.
\end{equation}
Under this constraint the following statement holds.\\

\noindent {\bf Theorem 2.} {\it For weights constrained by $p_1w_1+p_2w_2=W$ (or $(k_1/k)w_1+(k_2/k)w_2=W$ for a fixed weights distribution), $R_0^1$ and $R_0^2$ attain their maxima when $w_1=w_2=W$, and the maximum values for both is $\displaystyle{R_{0}=(k-1)r=\frac{(k-1)\tau W}{\tau W+\gamma}}$.}\\

\noindent The proof of this result is presented in Appendix C.

The above results suggest that for the same average link weight and when the one-to-one correspondence between $p_1$ and $k_1/k$, and $p_2$ and $k_2/k$ holds, the basic reproductive ratio is higher on networks with random weight distribution than on networks with a fixed weight distribution. This, however,  does not preclude the possibility of having a network with random weight distribution with smaller average weight exhibiting an $R_0$ value that it is bigger than the $R_0$ value corresponding to a network where weights are fixed and the average weight is higher. The direct implication is that it is not sufficient to know just the average link weight in order to draw conclusions about possible epidemic outbreaks on weighted networks; rather one has to know the precise weight distribution that provides a given average weight.

\begin{figure}
\centering
\epsfig{file=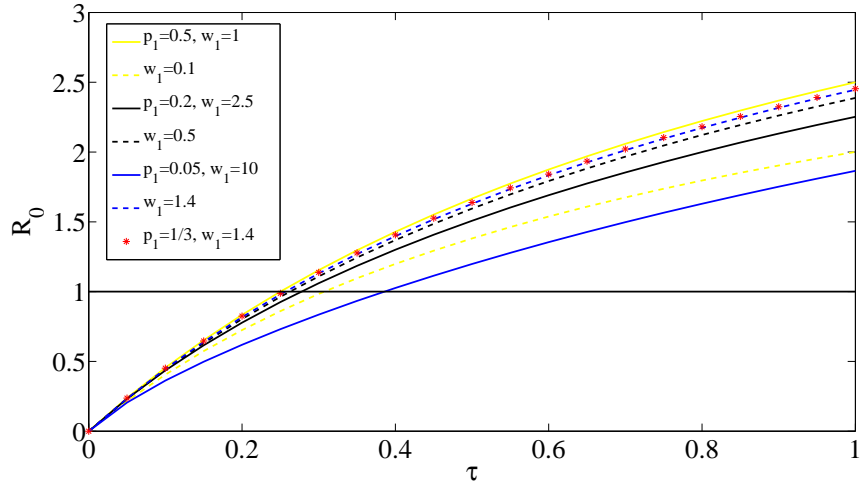,width=13cm}
\caption{Plots of the basic reproductive ratio $R_0$ for the two types of weighted networks with
different weight and weight frequency combinations, but with $p_1w_1+p_2w_2=1$. The case of homogenous networks with weights assigned at random (continuous lines with $p_1$ and $w_1$ given in the legend) considers the situation where the contribution of the two different weight types is equal ($p_1w_1=p_2w_2=0.5$) but with weight $w_1$ increasing (top to bottom) and its frequency decreasing. Increasing the magnitude of weights but reducing their frequency leads to smaller $R_0$ values. The case of homogeneous networks with fixed number of links of type $w_1$ and $w_2$ (dashed lines with only $w_1$ given in the legend) illustrates the situation where $w_1$ increases (bottom to top) while $p_1=k_1/k=1/3$ and $p_2=(k-k_1)/k=2/3$ remain fixed. Here the opposite tendency is observed with increasing weights leading to higher $R_0$ values. Finally, for the randomly distributed weights case, setting $p_1=1/3, w_1=1.4$ and observing $p_1w_1+p_2w_2=1$, we obtain $R_0$ ($\star$) values which compare almost directly to the fixed-weights case (top, dashed line). Other parameters are set to $k=6$, $k_1=2$ and $\gamma=1$.}
\label{rzero}
\end{figure}

Figure~\ref{rzero} shows how the basic reproductive ratio changes with the transmission rate $\tau$ for different weight distributions. When links on a homogeneous network are distributed at random, the increase in the magnitude of one specific link weight (e.g. $w_1$) accompanied by a decrease in its frequency leads to smaller $R_0$ values.  This is to be expected since the contribution of the different link types in this case is kept constant ($p_1w_1=p_2w_2=0.5$) and this implies that the overall weight of the network links accumulates in a small number of highly weighted links with most links displaying small weights and thus making transmission less likely. The statement above is more rigorously underpinned by the results of Theorem 1 \& 2 which clearly show that equal or more homogeneous weights lead to higher values of the basic reproductive ratio.  For the case of fixed weight distribution, the changes in the value of $R_0$ are investigated in terms of varying the weights, so that overall weight in the network remains constant. This is constrained by fixing values of $p_1$ and $p_2$ and, in this case, the highest values are obtained for higher values of $w_1$. The flexibility here is reduced due to $p_1$ and $p_2$ being fixed, and a different link breakdown may lead to different observations. The top continuous line in Fig.~\ref{rzero} corresponds to the maximum $R_0$ value achievable for both models if the $p_1w_1+p_2w_2=1$ constraint is fulfilled.

%%%%%%%
\subsection{$R_0$-like threshold for the $SIR$ model - a pairwise model perspective}
%%%%%%%
%A simple linear stability analysis of the disease-free steady state in terms of the pairwise model reveals an $R_0$ like threshold equal to $R_0=\frac{\tau (p_1w_1+p_2w_2)(k-1)}{\gamma}$ \textbf{(Is
%this the same if the second type of closure is used?)}.

To compute the value of $R_0$-like quantity from the  pairwise model, we use the approach suggested by Keeling \cite{Keeling1999}, which utilises the local spatial/network structure and correctly accounts for correlations between susceptible and infectious nodes early on in the epidemics. This can be achieved by looking at the early behaviour of $[SI]_1/[I]=\lambda_1$ and $[SI]_2/[I]=\lambda_2$ when considering links  of only two different weights. In line with Eames \cite{Eames2008}, we start from the evolution equation of $[I]$
%\begin{equation}
\[
\dot{[I]}=(\tau w_1[SI]_1/[I]+\tau w_2[SI]_2/[I] - \gamma)[I],
\]
%\end{equation}
where from the growth rate $\tau w_1\lambda_1+\tau w_2\lambda_2 - \gamma$ it is easy to define the threshold quantity $R$ as follows,
\begin{equation}\label{Rdef}
R=\frac{\tau w_1\lambda_1+\tau w_2\lambda_2}{\gamma}.
\end{equation}
For the classic closure (\ref{Clos_cl}), one can compute the early quasi-equilibria for $\lambda_1$ and $\lambda_2$ directly from the pairwise equations as follows
%\begin{equation}
\[
\lambda_{1} = \frac{\gamma{(k-1)}{p_{1}}{R}}{\tau{w_{1}}+\gamma{R}}\,\,\, \text{and} \,\,\,  \lambda_{2} = \frac{\gamma{(k-1)}(1-p_{1}){R}}{\tau{w_{2}}+\gamma{R}}.
\]
%\end{equation}
Substituting these into (\ref{Rdef}) and solving for $R$ yields
 \begin{equation}\label{Rtype1}
 R = \frac{R_{1}+R_{2}+\sqrt{(R_{1}+R_{2})^2+4R_{1}R_{2}Q}}{2},
 \end{equation}
 where
 \[
 \begin{array}{l}
 \displaystyle{R_{1} = \frac{\tau{w_{1}}[(k-1)p_{1}-1]}{\gamma},\hspace{0.3cm}
 R_{2} = \frac{\tau{w_{2}}[(k-1)p_{2}-1]}{\gamma},}\\\\%\hspace{0.3cm}
 \displaystyle{Q = \frac{k-2}{[(k-1)p_1-1][(k-1)p_2-1]},}
 \end{array}
 \]
with details of all calculations presented in Appendix D. We note that $R>1$ will result in an epidemic, while $R<1$ will lead to the extinction of the disease. It is straightforward to show that for equal weights, say $W$, the expression above reduces to $R=\tau W (k-2)/\gamma$ which is in line with $R_{0}$ value in \cite{Keeling1999} for unclustered, homogeneous networks. Under the assumption of a fixed total weight $W$, one can show that similarly to the network-based basic reproductive ratio, $R$ achieves its maximum when $w_1=w_2=W$.
%\textbf{(I believe that under the constant average weight constraint this will also be maximised when the weights are equal! Shall we test this? We could then say that this confirms that equal weights maximises $R$ and the bigger $R$ the more likely it is to observe an epidemic.)}.

In a similar way, for the modified closure (\ref{Clos}), we can use the same methodology to derive the threshold quantity as
\begin{equation}\label{Rtype2}
 R=\frac{R_{1}+R_{2}+\sqrt{(R_{1}+R_{2})^2+4R_{1}R_{2}(Q-1)}}{2},
\end{equation}
where
\[
 \displaystyle{R_{1} = \frac{\tau{w_{1}}(k_{1}-2)}{\gamma},\hspace{0.3cm}R_{2} = \frac{\tau{w_{2}}(k_{2}-2)}{\gamma},\hspace{0.3cm} Q = \frac{k_{1} k_{2}}{(k_{1}-2)(k_{2}-2)}.}
\]
For this closure once again, $R>1$ results in an epidemic, while for $R<1$, the disease dies out. Details of this calculations are shown in Appendix D. It is noteworthy that one can derive expressions (\ref{Rtype1}) and (\ref{Rtype2}) by considering the leading eigenvalue of linearization of system (\ref{SIR}) near its disease-free steady state with the corresponding pairwise closures given in (\ref{Clos_cl}) and (\ref{Clos}).

Finally, we note that this seemingly $R_0$-lookalike, $R=\tau W (k-2)/\gamma$ for the equal weights case $w_1=w_2=W$ is a multiple of $(k-2)$ as opposed to $(k-1)$ as is the case for the $R_0$ derived based on the individual-based perspective, where, for equal weights, $R_0^1=R_0^2=\tau W (k-1)/(\tau W + \gamma)$. This obviously highlights the strong dependency of $R_0$ on the modelling framework and also the difficulty in trying to reconcile findings based on different models.

%%%%%%%
\subsection{The performance of pairwise models and the impact of weight distributions on the dynamics of epidemics}
%%%%%%%

\begin{figure}

\hspace{-0.5cm}
\epsfig{file=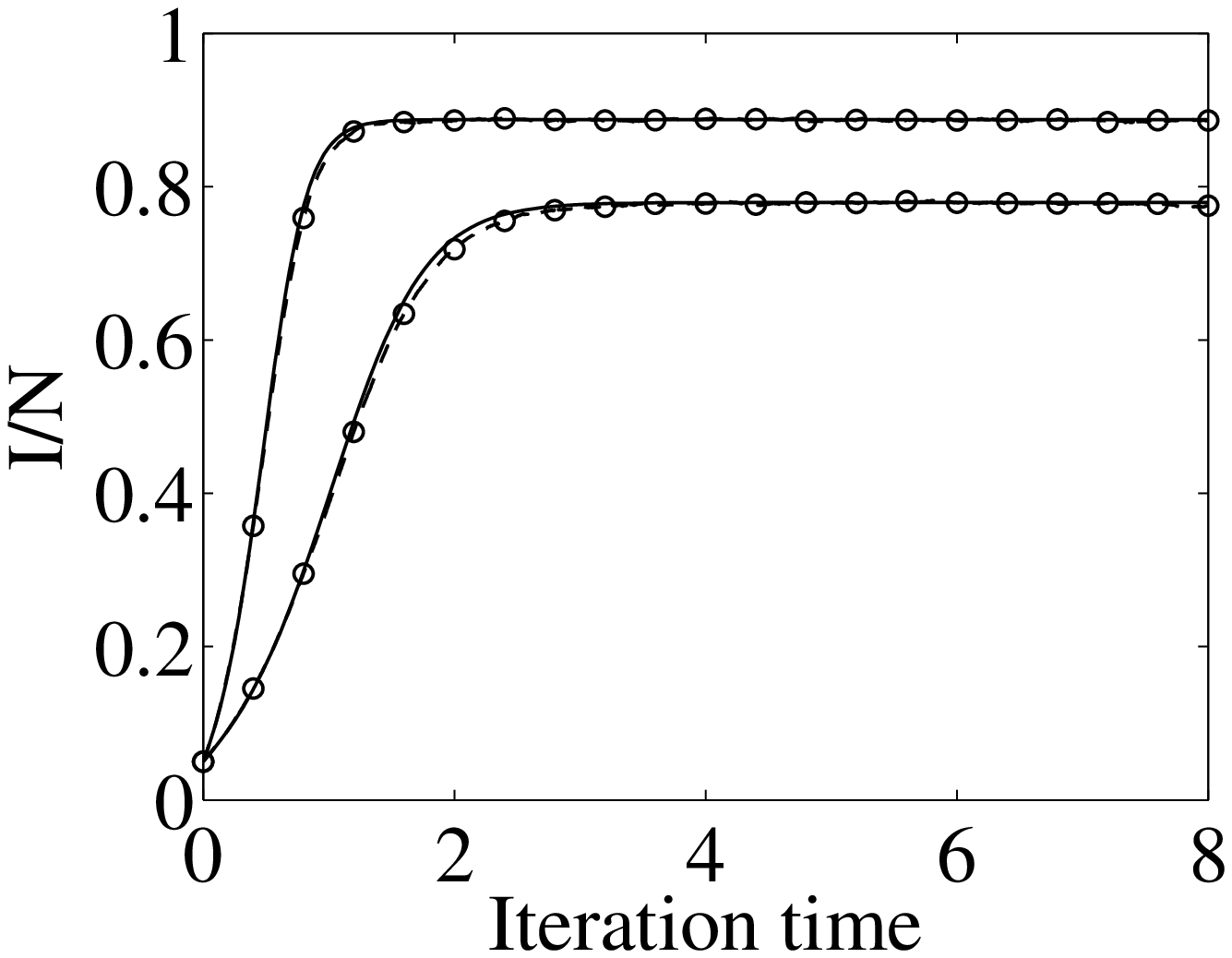,width=6.5cm}
%\hspace{0.5cm}
\epsfig{file=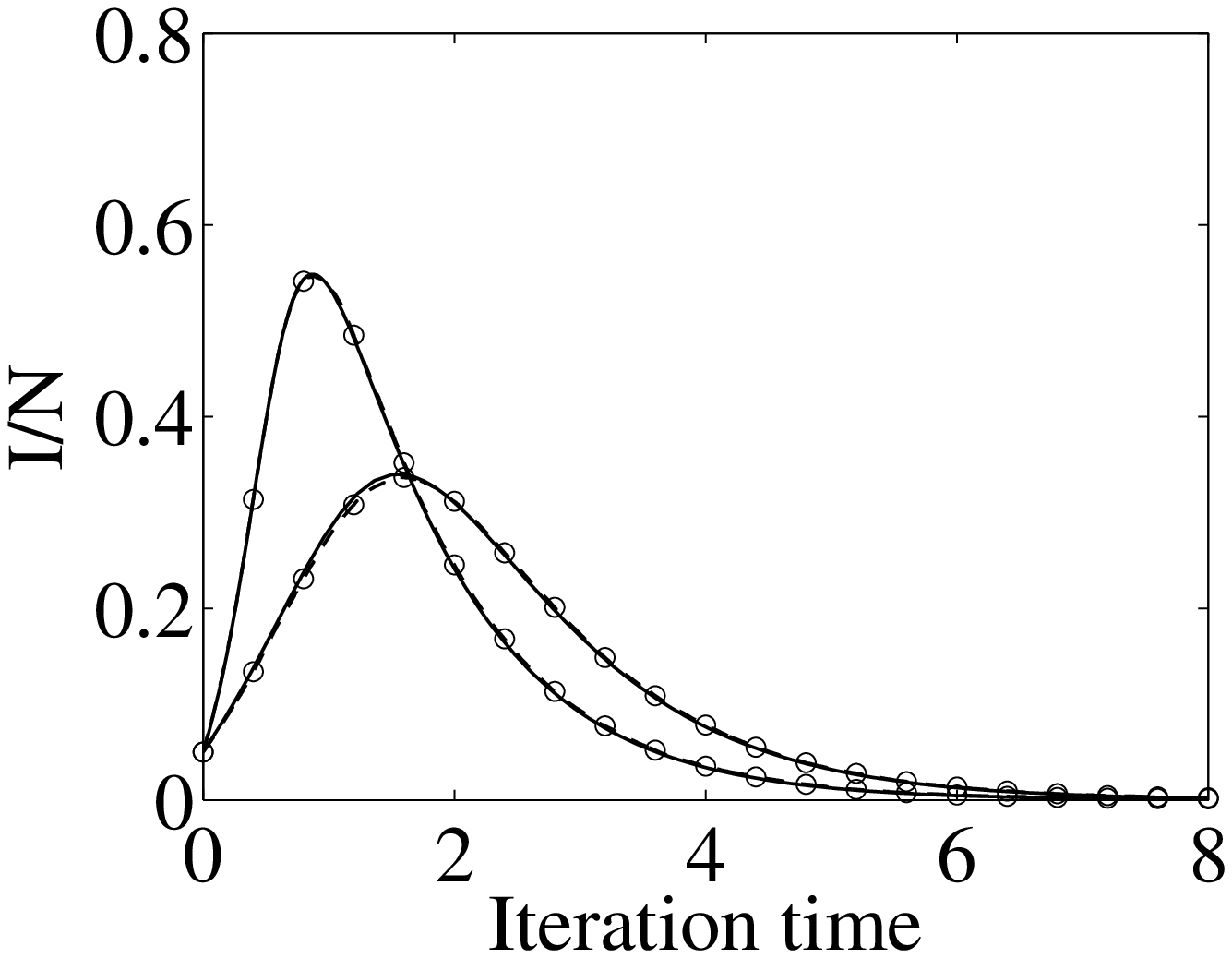,width=6.5cm}

%\caption{

%    \centering
%        \subfigure[$SIS$]
%    {
%        \includegraphics[width=7cm]{SIS_1a_Gilles_10ng_10run_new.eps}
%        \label{Fig2a}
%    }
%    \subfigure[$SIR$]
%    {
%        \includegraphics[width=7cm]{SIR_1b_Gilles_10ng_10run_new.eps}
%        \label{Fig2b}
%    }
    \caption{The infection prevalence ($I/N$) from the pairwise and simulation
models for homogeneous random networks with random weight distribution (ODE: solid line, simulation: dashed line and (o)). All nodes have degree $k = 5$ with
$N=1000$, $I_0=0.05N$, $\gamma = 1$ and $\tau = 1$. From top to bottom, the parameter values are: $w_{1}=5, p_1=0.2, w_{2}=1.25, p_2=0.8 $ (top), and $w_{1}=0.5, p_1=0.5, w_{2}
=1.5, p_2=0.5 $ (bottom). The left and right panels represent the $SIS$ and $SIR$ dynamics, respectively.}
\vspace{2cm}
\label{Fig2}
\end{figure}

To evaluate the efficiency of the pairwise approximation models, we will now compare numerical solutions of models (\ref{SIS}) and (\ref{SIR}) to
results obtained from the corresponding network simulation.
%We analyse in parallel the performance of the pairwise model for the $SIS$ and $SIR$ dynamics.
%This is achieved by a systematic analysis of the agreement between results from the pairwise and simulation model.
The discussion around the comparison of the two models is interlinked with the discussion of the impact of different weight distributions/patterns on the overall epidemic dynamics.
We being our numerical investigation by considering weight distributions with moderate heterogeneity. This is illustrated in Fig.~\ref{Fig2}, where excellent agreement between simulation and pairwise models is obtained. The agreement remains valid for both $SIS$ and $SIR$ dynamics, and networks with higher average link weight lead to higher prevalence levels at equilibrium for $SIS$ and higher infectiousness peaks for $SIR$.

\begin{figure}
%    \centering
%        \subfigure[$SIS$]
%    {
%        \includegraphics[width=7.5cm]{SIS_2a_Gilles_10ng_10run_new.eps}
%        \label{Fig3a}
%    }
%    \subfigure[$SIR$]
%    {
%        \includegraphics[width=7.5cm]{SIR_2b_Gilles_10ng_10run_new.eps}
%        \label{Fig3b}
%    }
\epsfig{file=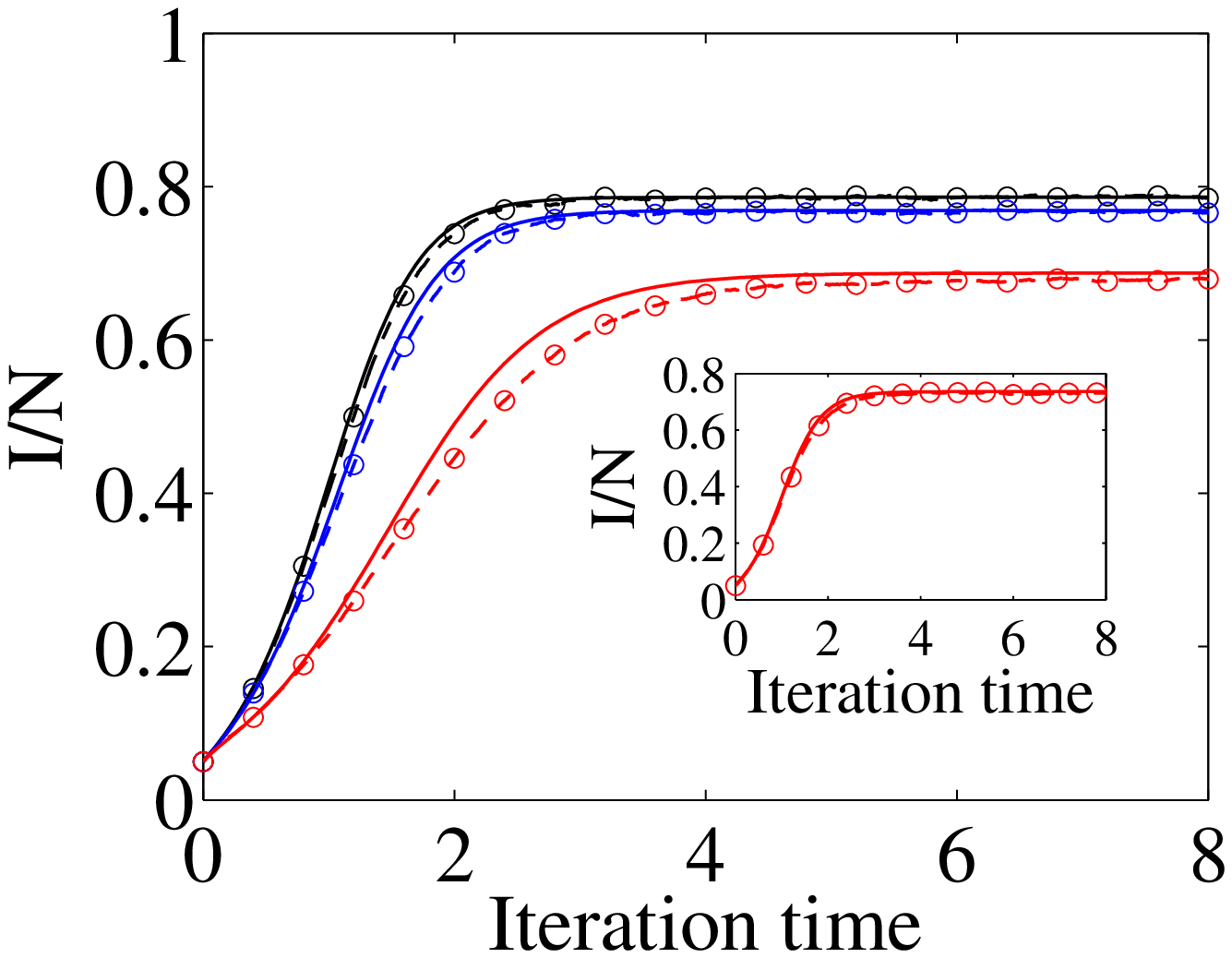,width=6.3cm}
%\hspace{0.5cm}
\epsfig{file=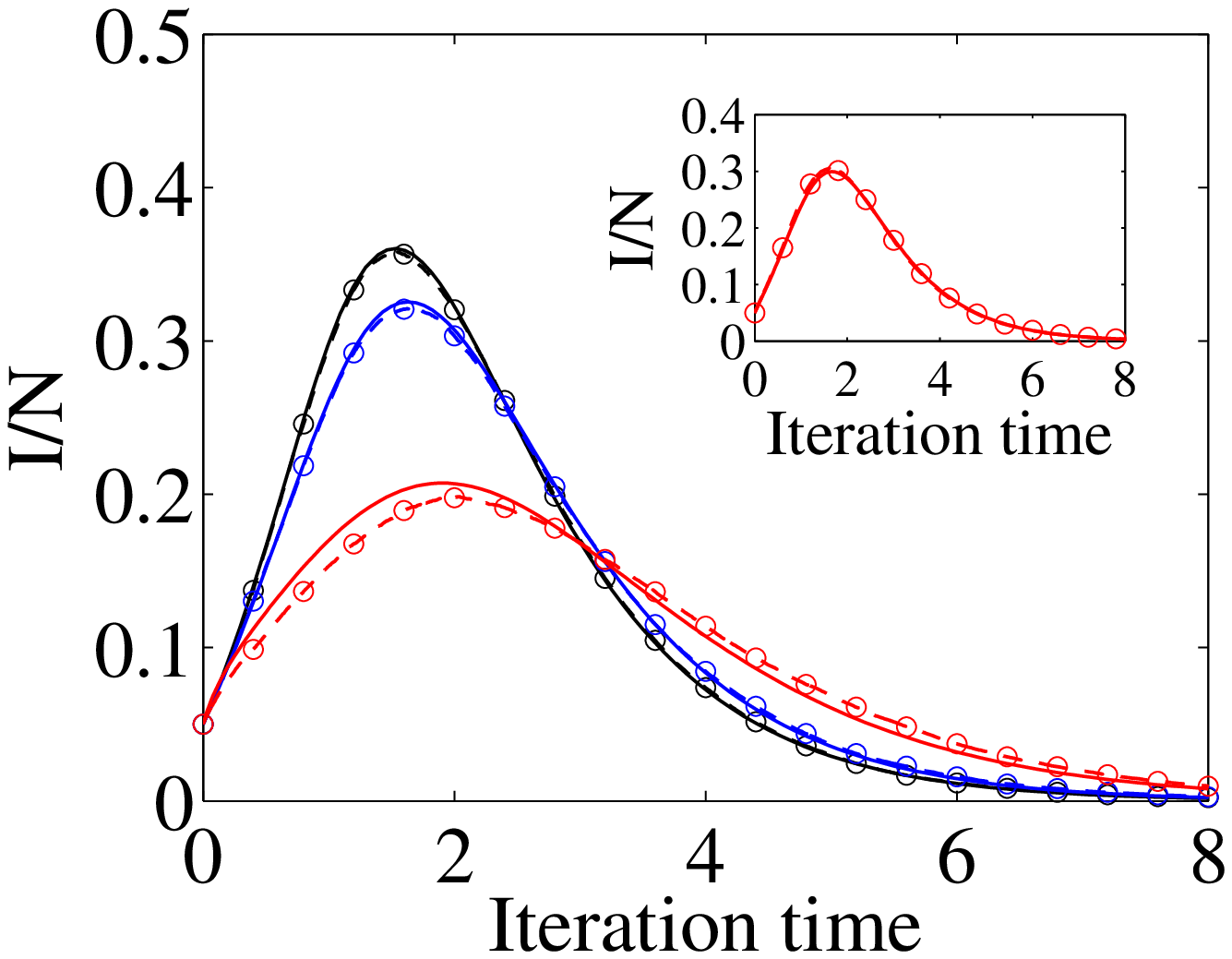,width=6.3cm}

    \caption{The infection prevalence ($I/N$) from the pairwise and simulation models for homogenous networks with random weight distribution (ODE: solid
line, simulation: dashed line and (o)). %which have the same average
%weight but difference the distribution weights.
All numerical tests use $N=1000$, $I_0=0.05N$, $k = 5$, $\gamma = 1, \tau = 1$ and $p_1=0.05$ ($p_2=1-p_1=0.95$).
From top to bottom, $w_{1} = 2.5,  5, 10$, $w_{2} = 0.875/0.95, 0.75/0.95, 0.5/0.95$.
The weight distributions are chosen such that the average link weight, $p_1w_1+p_2w_2=1$, remains constant.
Insets of (a) and (b): the same parameter values as for the lowest prevalence plots but,
with $k = 10$ and $\tau = 0.5$. The left and right panel represent the $SIS$ and $SIR$ dynamics, respectively.}
\label{Fig3}
\end{figure}

Next, we explore the impact of weight distribution under the
condition that the average weight remains constant (i.e.
$p_1w_1+p_2w_2=1$, where without loss of generality the average
weight has been chosen to be equal to 1). First, we keep the
proportion of edges of type one (i.e. with weight $w_1$) fixed and
change the weight itself by gradually increasing its magnitude.
Due to the constraint on the average weight and the condition
$p_2=1-p_1$, the other descriptors of the weight distribution
follow. Fig.~\ref{Fig3} shows that concentrating a large portion
of the total weight on a few links leads to smaller epidemics,
since the majority of links are low-weight and thus have a small
potential to transmit the disease. This effect is exacerbated for
the highest value of $w_1$; in this case $95\%$ of the links are
of weight $w_2=(1-p_1w_1)/(1-p_1)=0.5/0.95$ leading to epidemics
of smallest impact (Fig.~\ref{Fig3}(a)) and smallest size of
outbreak (Fig.~\ref{Fig3}(b)).

\begin{figure}
%    \centering
%    \subfigure[$SIS$]
%    {
%        \includegraphics[width=7cm]{SIS_3a_Gilles_10ng_10run_new.eps}
%        \label{Fig4a}
%    }
%    \subfigure[$SIR$]
%    {
%        \includegraphics[width=7cm]{SIR_3b_Gilles_10ng_10run_new.eps}
%        \label{Fig4b}
%    }
\epsfig{file=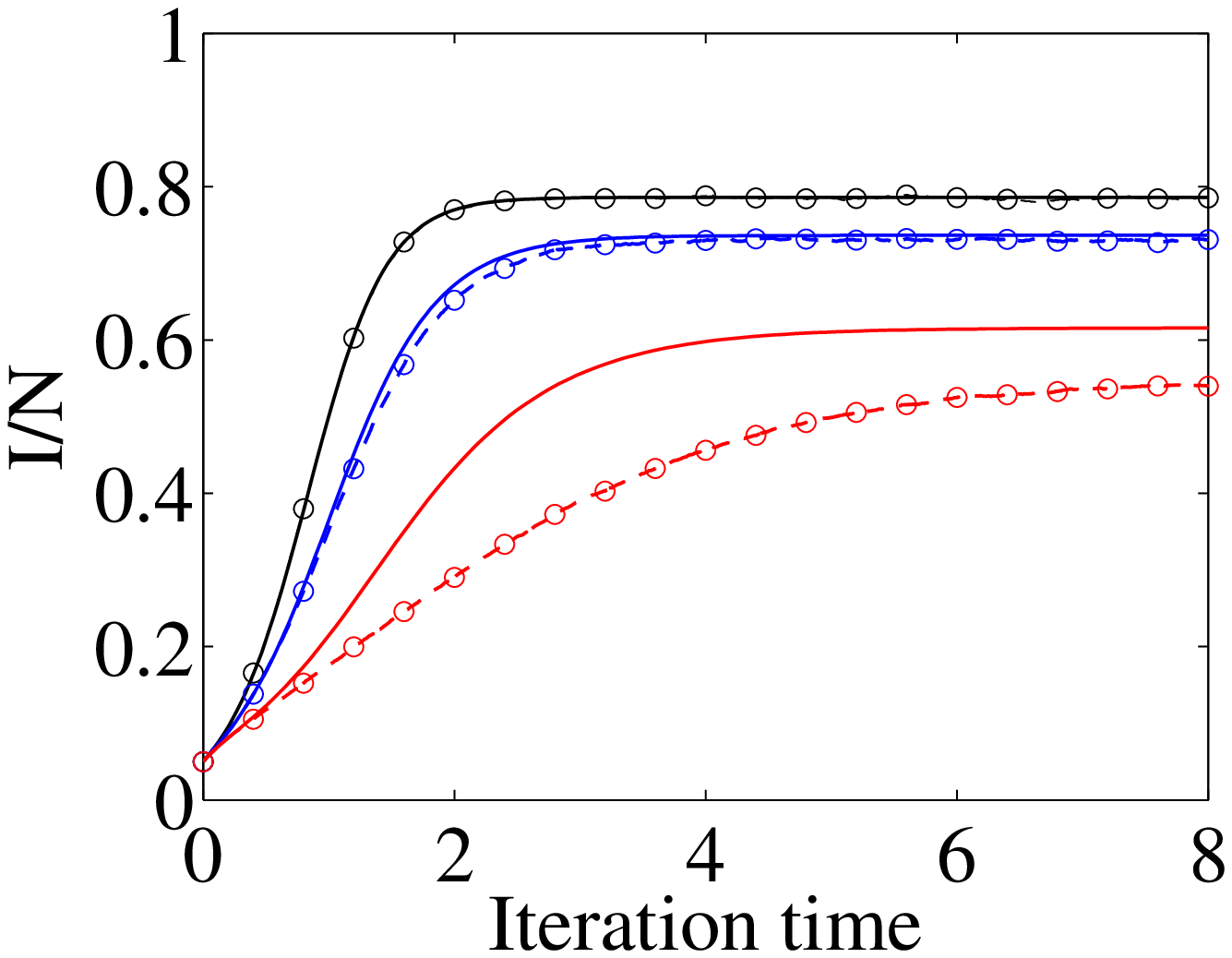,width=6.3cm}
%\hspace{0.5cm}
\epsfig{file=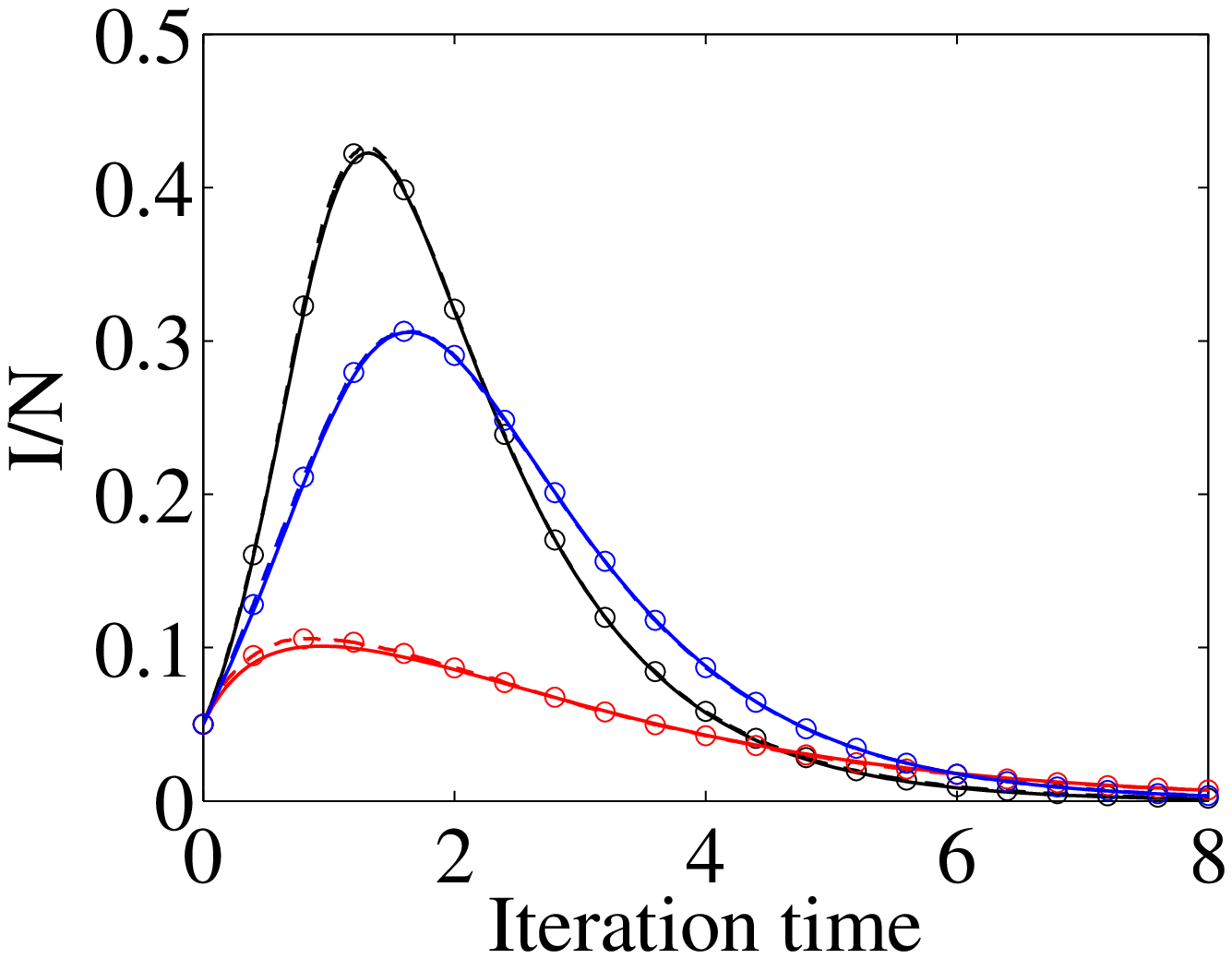,width=6.3cm}
    \caption{The infection prevalence ($I/N$) from the pairwise and simulation model for homogenous networks with random weight distribution (ODE: solid
line, simulation: dashed line and (o)).
% which have the same average
% weight but difference the probability of weights.
All numerical tests use $N=1000$, $I_0=0.05N$, $ k = 10$, $\gamma = 1, \tau =
0.5$ and $w_{1} = 10$. From top to bottom, $P(w_{1}) = 0.01, 0.05,
0.09$, $w_{2} = 0.9/0.99, 0.5/0.95, 0.1/0.91$.  Here also $p_2=1-p_1$ and $p_1w_1+p_2w_2=1$. The left and right panel represent the $SIS$ and $SIR$ dynamics, respectively.}
\vspace{2cm}
\label{Fig4}
\end{figure}

While the previous setup kept the frequency of links constant while changing the weights, one can also investigate the impact of keeping at least one of the weights constant (e.g. the larger one) and changing its frequency. To ensure a fair comparison, here we also require that the average link weight over the whole network is kept constant. When such highly weighted links are rare, the system approaches the non-weighted network limit where the transmission rate is simply scaled by $w_2$ (the most abundant link type). As Fig.~\ref{Fig4} shows, in this case, the agreement is excellent, and as the frequency of the highly weighted edges/links increases, disease transmission is less  severe.

Regarding the comparison of the pairwise and simulation models, we note that while the agreement is generally good for a large part of the disease and weight parameter space, the more extreme scenarios of weight distribution result in poorer agreement. This is illustrated in both Figs.~\ref{Fig3} and \ref{Fig4} (see bottom curves), with the worst agreement for the $SIS$ dynamics. The insets in Fig.~\ref{Fig3} show that increasing the average connectivity improves the agreement. However, the cause of disagreement is due to a more subtle effect driven also by the weight distribution. For example, in Fig.~\ref{Fig4}, the average degree in the network is $10$, higher then used previously and equal to that in the insets from Fig.~\ref{Fig3}, but despite this, the agreement is still poor.

\begin{figure}
%    \centering
%    \subfigure[$SIS$]
%    {
%        \includegraphics[width=7cm]{SIS_4a_Gilles_10ng_10run_new-1.eps}
%        \label{Fig5a}
%    }
%    \subfigure[$SIR$]
%    {
%        \includegraphics[width=7cm]{SIR_4b_Gilles_10ng_10run_new-1.eps}
%        \label{Fig5b}
%    }
\epsfig{file=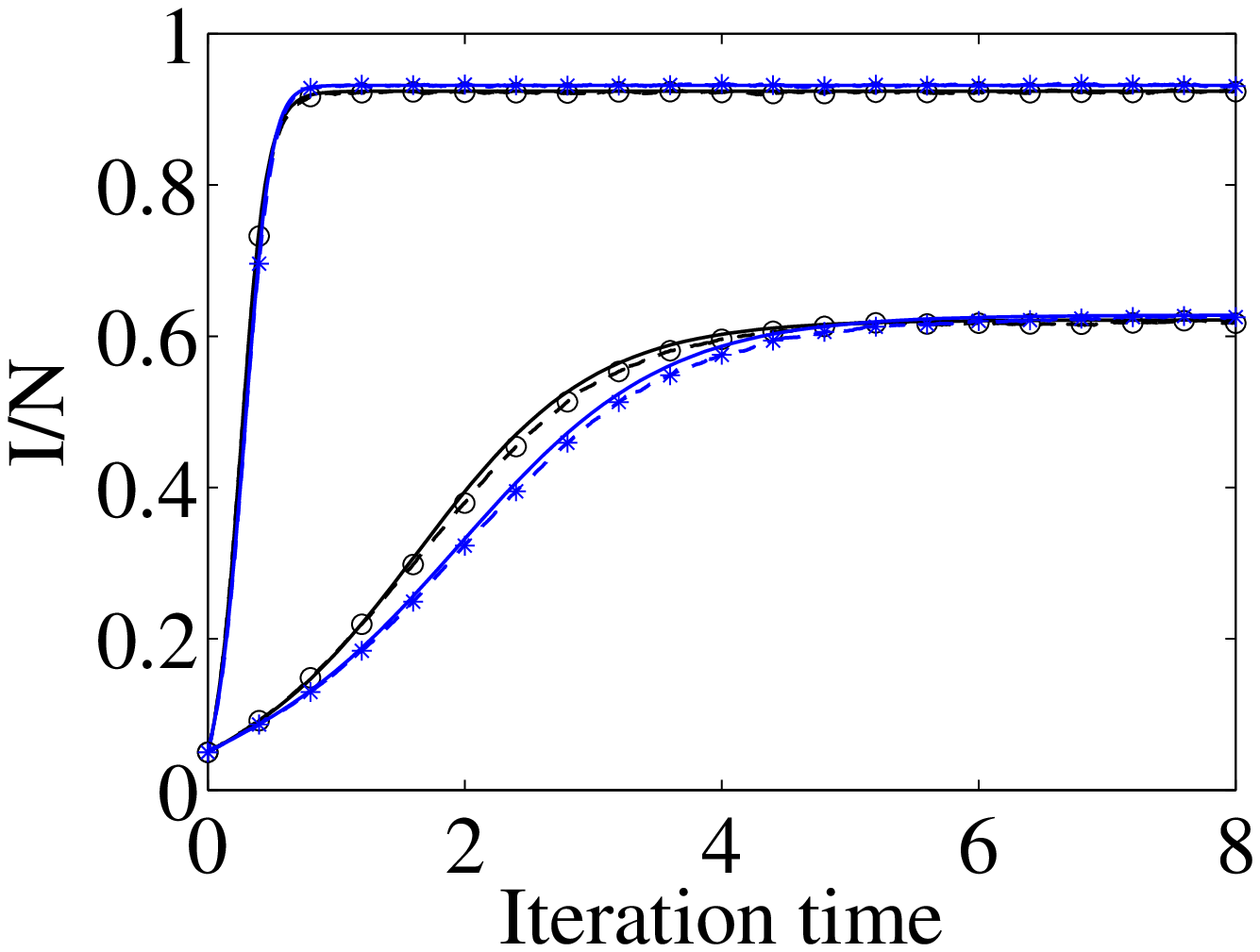,width=6.3cm}
%\hspace{0.5cm}
\epsfig{file=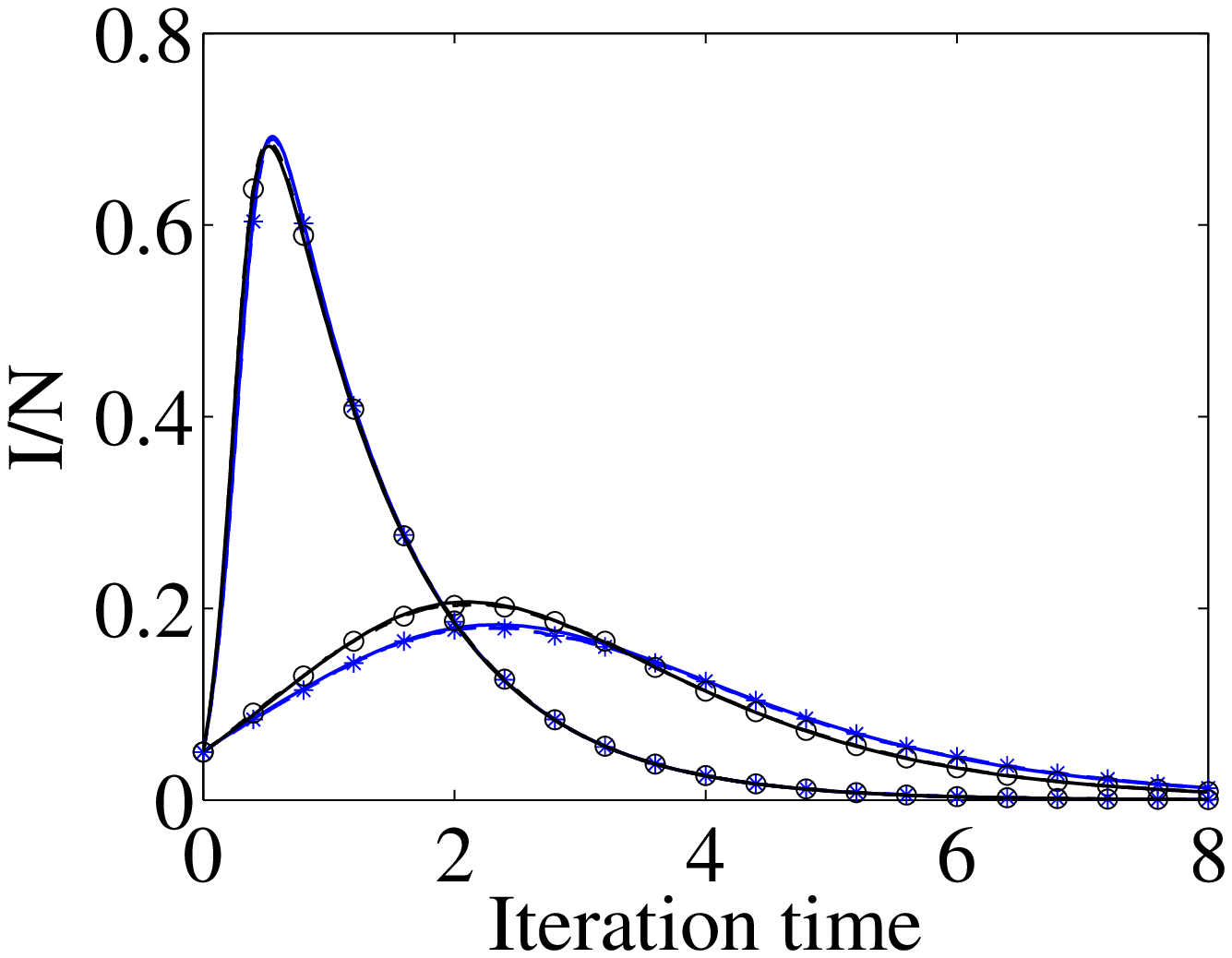,width=6.3cm}
\caption{The infection prevalence ($I/N$) based on random (model 1) and fixed (model 2) weight distribution (ODE: black (1) and blue (2) solid line, simulation results: same as ODE but dashed lines,
and ($\circ$) and ($\ast$)). All numerical tests use $N=1000$, $I_0=0.05N$, $k = 10$, $k_1=2$, $k_2=8$, $p_1=k_1/k$, $p_2=k_2/k$, $w_1=10$, $w_2=1.25$ and $\gamma = 1$. The rate of infection $\tau=0.5$ (top) and $\tau=0.1$ (bottom). The left and right panel represent the $SIS$ and $SIR$ dynamics, respectively.}
\label{Fig5}
\end{figure}

The two different weighted network models are compared in Fig.~\ref{Fig5}. This is done be using the same link weights and setting $p_1=k_1/k$ and $p_2=k_2/k$.
Epidemics on network with random weight distribution grow faster and, given the same time scales of the epidemic, this is line with results
derived in Theorem 1 \& 2 and findings concerning the growth rates. The difference is less marked for larger values of $\tau$ where a significant proportion of the
nodes becomes infected.

\begin{figure}
%    \centering
%    \subfigure[$SIS$]
%    {
%        \includegraphics[width=7cm]{SIS_5a_Gilles_10ng_10run_new.eps}
%        \label{$SIS$Fig6a}
%    }
%    \subfigure[$SIR$]
%    {
%        \includegraphics[width=7cm]{SIR_5b_Gilles_10ng_10run_new.eps}
%        \label{Fig6b}
%    }
\epsfig{file=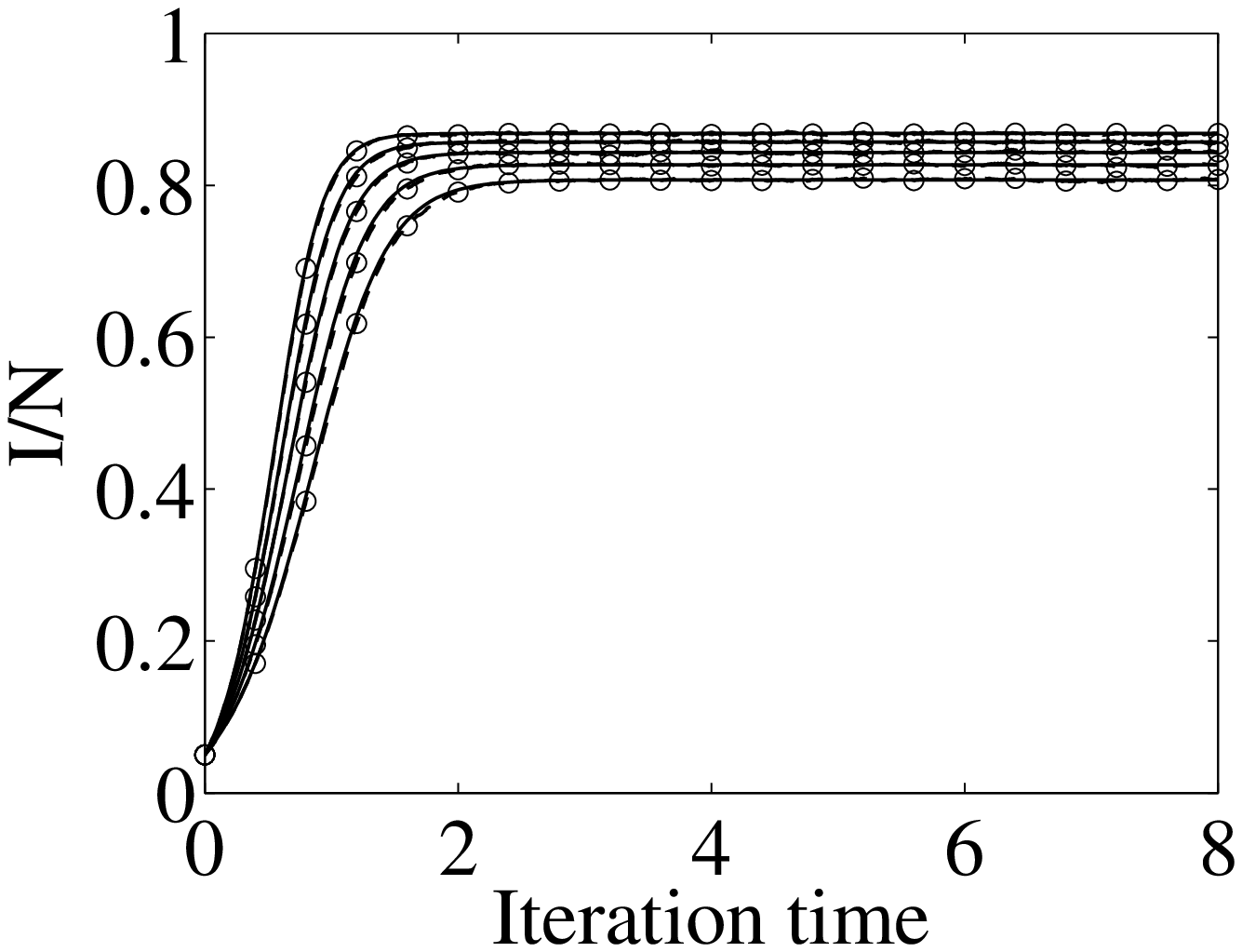,width=6.3cm}
%\hspace{0.5cm}
\epsfig{file=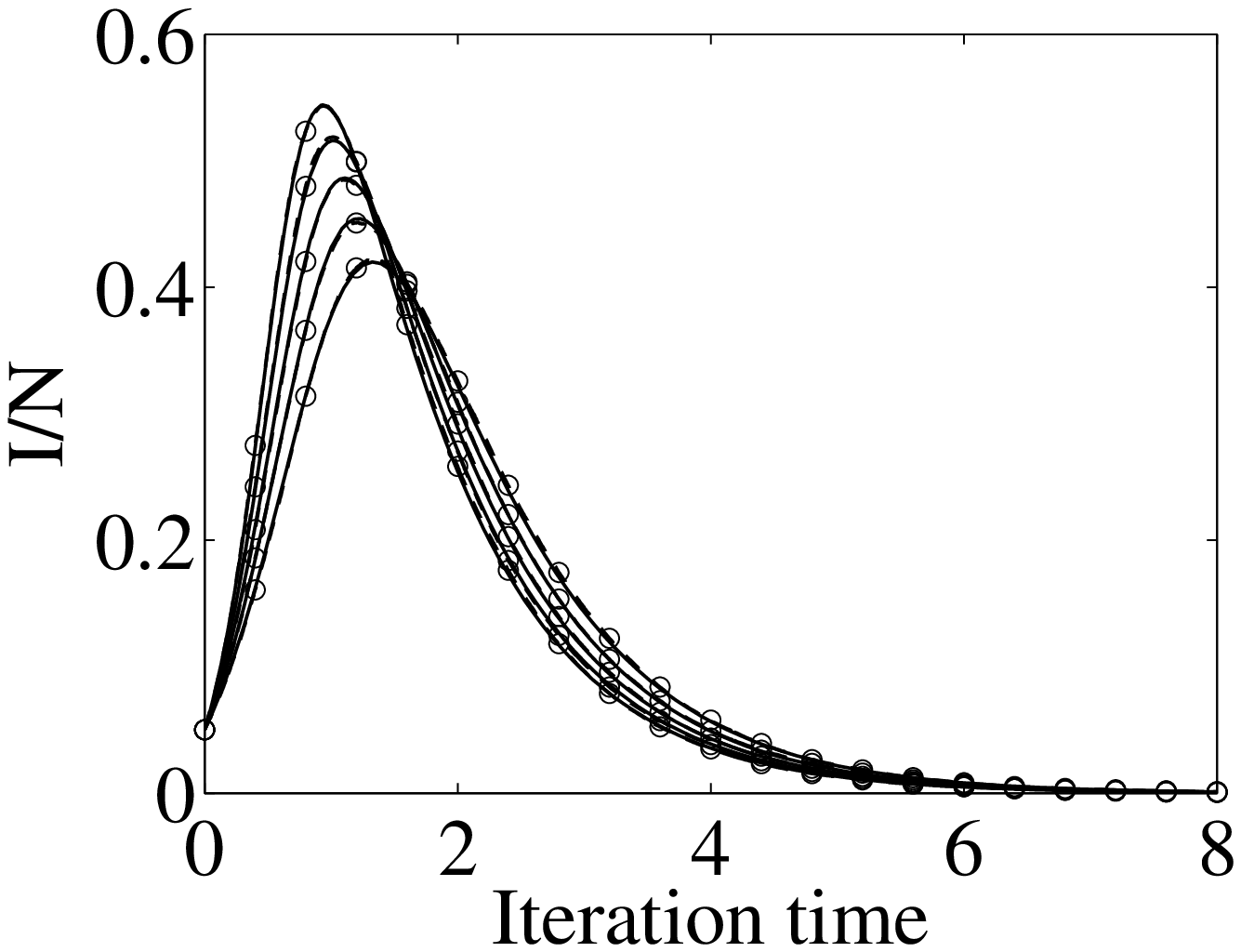,width=6.3cm}
    \caption{The infection prevalence ($I/N$) for a fixed weight distribution (ODE: solid lines, simulation results: dashed lines and (o)). All numerical tests use $N=1000$, $I_0=0.05N$, $ k = 6$,
$\gamma = 1, \tau = 1$ and $w_{1} = 1.4, w_{2} = 0.8$. From top
to bottom : $k_{1} = 5, 4, 3, 2, 1$ and $k_{2} = k-k_{1}$. The left and right panel represent the $SIS$ and $SIR$ dynamics, respectively.}
\vspace{2cm}
\label{Fig6}
\end{figure}

In Fig.~\ref{Fig6} the link weight composition is altered by decreasing the proportion of highly-weighted links. As expected, the reduced average link weight across the network
leads to epidemics of smaller size while keeping the excellent agreement between simulation and pairwise model results.

%\textbf{Por, please follow what I have written above and write a discussion for figures 5 and 6!!!}

%%%%%%%%%%%%%%%%%%%%%%%%%%%%%%%%%%%%%%%%%%%%%%%%%%%%%%%%%%%%%%%%%%%%%%%%%%%%%%%%%%%%%%%%%%%%%%%%%%
\section{Discussion}
%%%%%%%%%%%%%%%%%%%%%%%%%%%%%%%%%%%%%%%%%%%%%%%%%%%%%%%%%%%%%%%%%%%%%%%%%%%%%%%%%%%%%%%%%%%%%%%%%%

%\textbf{The useful thing that comes out of this, I think, is that the added heterogeneity of link weights doesn't behave in the same way as most other heterogeneities in models: usually, heterogeneities
%lead R to increase but potentially for final size to fall. Here, they concentrate infectiousness on fewer target individuals leading to a fall in $R_0$ (though it seems model 2 goes the other way - will have
%to think about that!). The higher the weight, the more "networky" the network, so it's taking it further from the mass-action type of model and lowering $R_0$ as it does so. Nice.}

The present study has explored the impact of weight heterogeneity and highlighted that the added heterogeneity of link weights does not manifest itself in the same way as most other heterogeneities in epidemic models on networks. Usually, heterogeneities lead to an increase in $R_0$ but potentially for final size to fall. However, for weighted networks the concentration of infectiousness on fewer target link, and thus target individuals, leads to a fall in $R_0$ for both homogeneous random and fixed weight distribution models. Increased heterogeneity in weights accentuates the locality of contact and is taking the model further from the mass-action type models. Infection is concentrated along a smaller number of links, which results in wasted infectivity and lower $R_0$. This is in line with similar results \cite{BDL,r13,Yan} where different modelling approaches have been used to capture epidemics on weighted networks.

The models proposed in this paper are simple mechanistic models with basic weight distributions, but despite this they provide a good basis for analysing disease dynamics on weighted networks in a rigorous and systematic way. The modified pairwise models have performed well, and provide good approximation to direct simulation. As expected, the agreement with simulations typically breaks down at or close to the threshold, but away from it, pairwise models provide a good counterpart or alternative to simulation. Disagreement only appears for extreme weight distributions, and we hypothesise that this is mainly due to the network becoming more modular with islands of nodes connected by links of low weight being bridged together by highly weighted links. A good analogy to this is provided by considering the case of a pairwise model on unweighted networks specified in terms of two network metrics, node number $N$ and average number of links $k$. The validity of the pairwise model relies on the network being connected up at random, or according to the configuration model. This can be easily broken by creating two sub-networks of equal size both exhibiting the same average connectivity. Simulations on such type of networks will not agree with the pairwise model, and highlights that the network generating algorithm can push the network out of the set of `acceptable' networks. We expect that this or similar argument can more precisely explain why the agreement breaks down for significant link-weight heterogeneity.

\begin{figure}
%\begin{center}
%\includegraphics[width=9cm]{Tau_and_I.eps}
\hspace{1cm}
\epsfig{file=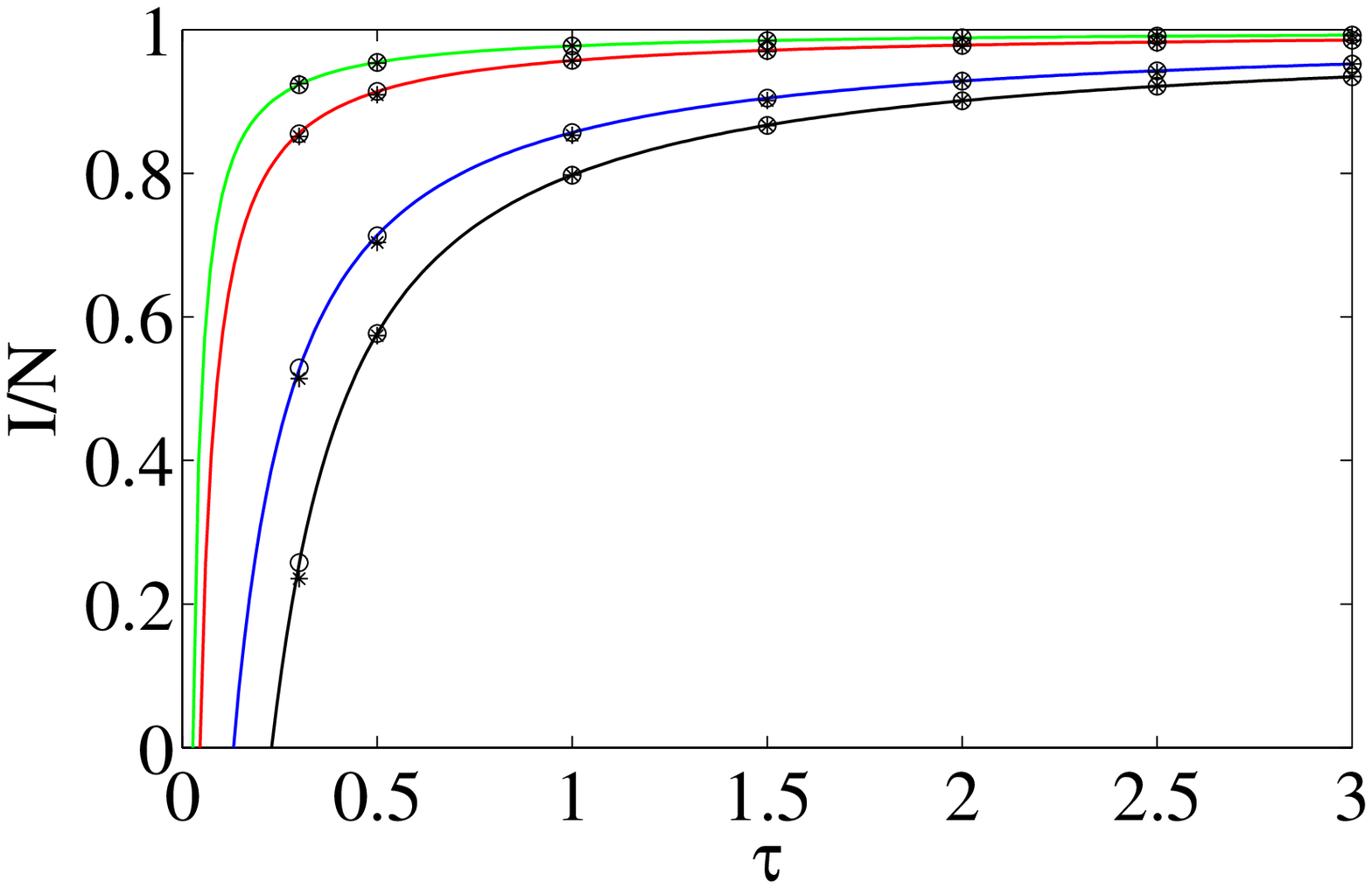,width=9cm}
\caption{Endemic steady state from the $SIS$ model on networks with random weight distribution.
The continuous lines correspond to the steady state computed numerically by setting all evolution equations in the pairwise system to zero.
These are complemented by finding the endemic steady state through direct integration of the ODE system for a long-enough time ($\circ$), as well as direct simulation ($\ast$). The first marker corresponds
to  $\tau = 0.3$ followed by $\tau=0.5, 1.0, \dots, 3.0$. All results are based on: $ k = 5$, $\gamma = 1$ and
$w_{1} = 10, w_{2} = 1$. From top to bottom : $p_{1} =0.9, 0.5,
0.1, 0.01$ and $p_{2} = 1-p_{1}$.}
\label{Kostia}
%\end{center}
\end{figure}

The usefulness of pairwise models is illustrated in Fig.~\ref{Kostia}, where the $I/N$ values are plotted for a range of $\tau$ values and for different weight distributions.
Here, the equilibrium value has been computed by finding the steady state directly from the ODEs (\ref{SIS}) by finding numerically the steady state solution of a set on nonlinear equations (i.e. $\dot{[A]}=0$ and $\dot{[AB]}=0$). To test the validity, the long term solution of the ODE is plotted along with results based on simulation. The agreement away from the threshold is excellent and illustrates clearly the impact of different weight distributions on the magnitude of the endemic threshold.

The models proposed here can be extended in a number of different ways. One potential avenue for further research is the analysis of correlations between link weight and node degree. This direction has been explored but in the context of classic compartmental mean-field models based on node degree \cite{NetwWithSaturation,LewiStone}. Given that pairwise models extend to heterogeneous networks such avenues can be further explored to include different type of correlations or other network dependent weight distributions. While this is a viable direction, it is expected that the extra complexity will make the pairwise models more difficult to analyse and disagreement between pairwise and simulation model more likely. Another theoretically interesting and practically important aspect is the consideration of different types of time delays, representing latency or temporary immunity \cite{BK10}, and the analysis of their effects on the dynamics of epidemics on weighted networks. The methodology presented in this paper can be of wider relevance to studies of other natural phenomena where overlay networks provide effective description. Examples of such systems include the simultaneous spread of two different diseases in the same population \cite{BK05}, the spread of the same disease but via different routes \cite{KissMultipRoutes} or the spread of epidemics concurrently with information about the disease \cite{Vasilis,KissInfo}. These areas offer other important avenues for further extensions.

%%%%%%%%%%%%%%%
\section*{Acknowledgements}
%%%%%%%%%%%%%%%
P. Rattana acknowledges funding for her PhD studies from the Ministry of Science and Technology, Thailand.
\newpage

%%%%%%%%%%
\section{Appendix}
%%%%%%%%%

%\textbf{Not convinced this adds anything. What are we trying to say?}

%%%%%%%%%%%%%%%%%%%%%%%%%%%%%%%%%%%%%%%%%%%%%%%
\subsection{Appendix A - Reducing the weighted pairwise models to the unweighted equivalents}\label{AppA}
%%%%%%%%%%%%%%%%%%%%%%%%%%%%%%%%%%%%%%%%%%%%%%%
We start from the system
\begin{equation}\label{edt1}
\begin{array}{l}
   \dot{[S]}= \gamma{[I]}-\tau{\sum_{n=1}^M {w_{n}}{[SI]_{n}}},\\\\
   \dot{[I]}= \tau{\sum_{n=1}^M {w_{n}}{[SI]_{n}}}-\gamma{[I]},\\\\
\ [\dot{SI}]_{m}= \gamma{([II]_{m}-[SI]_{m})}+\tau{\sum_{n=1}^M {w_{n}}{([SSI]_{mn}-[ISI]_{nm})}}-\tau{w_{m}}{[SI]_{m}},\\\\
   \ [\dot{II}]_{m}= -2{\gamma{[II]_{m}}}+2{\tau}{\sum_{n=1}^M w_{n}{[ISI]_{nm}}}+2{\tau}{w_{m}}{[SI]_{m}},\\\\
   \ [\dot{S}{S}]_{m}= 2{\gamma{[SI]_{m}}}-2{\tau{\sum_{n=1}^M w_{n}{[SSI]_{mn}}}},
   \end{array}
 \end{equation}
where $m=1,2, \dots, M$. To close this system of equations at the level of
pairs, we use the approximations
%\begin{eqnarray}\label{edt2}
\[
   \ [ABC]_{mn}=\frac{k-1}{k}\frac{[AB]_{m}[BC]_{n}}{[B]}.
\]
%\end{eqnarray}
To reduce these equations to the standard pairwise model for unweighted networks we use the fact that $\sum_{m=1}^M {[AB]_{m}} = [AB]$ for $A,B \in \{S, I\}$
and aim to derive the evolution equation for $[AB]$. Assuming that all weights are equal to some $W$, the following relations hold,
\begin{align*}
\dot{[SI]}&=\sum_{m=1}^M\dot{[SI]_m}\\
              &=\sum_{m=1}^M \left(\gamma{([II]_{m}-[SI]_{m})}+\tau{\sum_{n=1}^M {w_{n}}{([SSI]_{mn}-[ISI]_{nm})}}-\tau{w_{m}}{[SI]_{m}}\right)\\
              &=\gamma([II]-[SI])-\tau W [SI] +\tau W \sum_{m=1}^M {\sum_{n=1}^M{([SSI]_{mn}-[ISI]_{nm})}},
\end{align*}
where the summations of the triples can be resolved as follows,
\begin{align*}
\sum_{m=1}^M \sum_{n=1}^M[SSI]_{mn}&=\frac{k-1}{k}\sum_{m=1}^M [SS]_{m}\sum_{n=1}^M \frac{[SI]_{n}}{[S]}\\\\
               &=\frac{k-1}{k}\frac{[SS][SI]}{[S]}=[SSI].
 \end{align*}
Using the same argument for all other triples, the pairwise model for weighted networks with all weights being equal (i.e. $W=1$) reduces to the classic pairwise model, that is
\[
\begin{array}{l}
   \dot{[S]}= \gamma{[I]}-\tau{[SI]},\\\\
   \dot{[I]}= \tau{[SI]}-\gamma{[I]},\\\\
   \ \sum_{m=1}^M {[\dot{S}{I}]}=[\dot{S}{I}]= \gamma{([II]-[SI])}+\tau{[SSI]-[ISI]-[SI]},\\\\
   \ \sum_{m=1}^M {[\dot{I}{I}]}=[\dot{I}{I}]= -2{\gamma{[II]}}+2{\tau{([ISI]+[SI])}},\\\\
   \ \sum_{m=1}^M {[\dot{S}{S}]}=[\dot{S}{S}]= 2{\gamma{[SI]}}-2{\tau{[SSI]}}.
   \end{array}
   \]
A similar argument holds for the pairwise model on weighted networks with $SIR$ dynamics.
%%%%%%%%%%%%%%%%%%%%%%%%%%%%%%
\subsection{Appendix B - Proof of Theorem 1}
%%%%%%%%%%%%%%%%%%%%%%%%%%%%%%
We illustrate the main steps needed to complete the proof of Theorem 1.
This revolves around starting from the inequality itself and showing via a series of algebraic manipulations that
it is equivalent to a simpler inequality that holds trivially. Upon using that $p_1k=k_1$, $p_2k=k_2$ and $p_2+p_1=1$, the original inequality can be rearranged to give
%\begin{eqnarray}
%\sqrt{((k_1-1)r_1-(k_2-1)r_2)^2+4k_1k_2r_1r_2} \le (k_1-1)r_1+(k_2-1)r_2+2r_1p_1+2r_2p_1
%\end{eqnarray}
\begin{equation}
\sqrt{[(k_1-1)r_1-(k_2-1)r_2]^2+4k_1k_2r_1r_2} \le (k_1-1)r_1+(k_2-1)r_2+2r_1p_2+2r_2p_1.
\end{equation}
Based on the assumptions of the Theorem, the right-hand side is positive, and thus this inequality is equivalent to
the one where both the left- and right-hand sides are squared. Combined with the fact that $p_2=1-p_1$,  after
a series of simplifications and factorizations this inequality can be recast as
\begin{equation}
4p_1(1-p_1)(r_1^2+r_2^2)+8kp_1(1-p_1)r_1r_2 \le 4kp_1(1-p_1)(r_1^2+r_2^2)+8p_1(1-p_1)r_1r_2,
\end{equation}
which can be further simplified to
\begin{equation}
4p_1(1-p_1)(r_1-r_2)^2(k-1) \ge 0,
\end{equation}
which holds trivially and thus completes the proof. We note that in the strictest mathematical sense the condition of the Theorem should be 
$(k_1-1)r_1+(k_2-1)r_2+2r_1p_2+2r_2p_1 \ge 0$. This holds if the current assumptions are observed since these are stronger but 
follow from a practical reasoning whereby for the network with fixed weight distribution, a node should have at least one link with every possible weight type.

%%%%%%%%%%%%%%%%%%%%%%%%%%%%%%
\subsection{Appendix C - Proof of Theorem 2}
%%%%%%%%%%%%%%%%%%%%%%%%%%%%%%
First, we show that $R_0^1$ is maximised when $w_1=w_2=W$. $R_0^1$ can be rewritten to give
 \begin{eqnarray}
       \ R_{0}^1 &=& (k-1)\left(p_{1}\frac{\tau{w_{1}}}{\tau{w_{1}}+r}+(1-p_{1})\frac{\tau{w_{2}}}{\tau{w_{2}}+r}\right).
 \end{eqnarray}
Maximising this given the constraint $w_{1}p_{1}+w_{2}(1-p_{1}) = W$ can be achieved by considering $R_0^1$ as a function of the two weights and incorporating the constraint into it
via the Lagrange multiplier method. Hence, we define a new function $f(w_1,w_2,\lambda)$ as follows
\[
\begin{array}{l}
\displaystyle{       \ f(w_{1},w_{2},\lambda) =
       (k-1)\left(p_{1}\frac{\tau{w_{1}}}{\tau{w_{1}}+r}+(1-p_{1})\frac{\tau{w_{2}}}{\tau{w_{2}}+r}\right)}\\\\
       \hspace{2cm}+\lambda(w_{1}p_{1}+w_{2}(1-p_{1})-W).
\end{array}
\]
Finding the extrema of this functions leads to a system of three equations
\[
\begin{array}{l}
  \displaystyle{     \ \frac{\partial{f}}{\partial{w_{1}}} = \frac{(k-1)p_{1}\tau\gamma}{(\tau{w_{1}}+\gamma)^2}+\lambda{p_{1}} =0, }\\\\
     \displaystyle{      \ \frac{\partial{f}}{\partial{w_{2}}} = \frac{(k-1)(1-p_{1})\tau\gamma}{(\tau{w_{2}}+\gamma)^2}+\lambda(1-p_{1}) =0, }\\\\
      \displaystyle{     \ w_{1}p_{1}+w_{2}(1-p_{1})-W = 0.}
\end{array}
\]
Expressing $\lambda$ from the first two equations and equating these two expressions yields
\begin{eqnarray}
       \ \frac{(k-1)\tau\gamma}{(\tau{w_{1}}+\gamma)^2} = \frac{(k-1)\tau\gamma}{(\tau{w_{2}}+\gamma)^2}.
\end{eqnarray}
Therefore,
\begin{eqnarray}
       \ w_{1} = w_{2} = W,
\end{eqnarray}
and it is straightforward to confirm that this is a maximum.

Performing the same analysis for $R_0^2$ is possible but it is more tedious. Instead, we propose a more elegant argument to show that $R_0^2$ under the constraint of constant total link weight achieves its maximum when $w_1=w_2=W$. The argument starts by considering $R_0^2$ when  $w_{1} = w_{2}=W$. In this case, and using that $r_2=r_1=r=\tau W/(\tau W + \gamma)$ we can write,
\[
\begin{array}{l}
 \displaystyle{        \ R^{2\ast}_0 = \frac{(k_1-1)r_1+(k_2-1)r_2+\sqrt{[(k_1-1)r_1-(k_2-1)r_2]^2+4k_1k_2r_1r_2}}{2}}\\\\
   \displaystyle{      \ = \frac{r (k_1+k_2-2)+\sqrt{r^2 [(k_1-1)-(k_2-1)]^2+4r^2 k_1 k_2}}{2}}\\\\
   \displaystyle{      \ = \frac{r (k_1+k_2-2)+r \sqrt{(k_1+k_2)^2}}{2}}\\\\
   \displaystyle{      \ = \frac{r (2k_1+2k_2-2)}{2} = r(k_1+k_2-1)=(k-1)r.}
\end{array}
\]
However, it is known from Theorem 1 that $R_0^2 \leq R_0^1$, and we have previously shown that $R_0^1$ under the present constraint achieves its maximum when $w_1=w_2=W$, and its maximum is equal to $(k-1)r$.
All the above can be written as
\begin{equation}
R_0^2 \leq R_{0}^1 \leq (k-1)r.
\end{equation}
Now taking into consideration that $R^{2\ast}_0=(k-1)r$, the inequality above can be written as
\begin{equation}
R_0^2 \leq R_{0}^1 \leq (k-1)r =R^{2\ast}_0,
\end{equation}
and this concludes the proof.

%%%%%%%%%%%%%%%%%%%%%%%%%%%%
\subsection{Appendix D - The $R_0$-like threshold $R$}
%%%%%%%%%%%%%%%%%%%%%%%%%%%%
Let us start from the evolution equation for $[I](t)$,
\begin{align*}
\dot{[I]}&= \tau{(w_{1}{{[SI]_{1}}}+w_{2}{{[SI]_{2}}})}-\gamma{[I]}\\
            &= \left[\tau{w_{1}}\left(\frac{[SI]_{1}}{[I]}\right)+\tau{w_{2}}\left(\frac{[SI]_{2}}{[I]}\right)-\gamma \right]{[I]}\\
            &= (\tau{w_{1}}\lambda_1+\tau{w_{2}}\lambda_2-\gamma){[I]},
 \end{align*}
where $ \lambda_{1}= \frac{[SI]_{1}}{[I]}$ and $ \lambda_{2}=\frac{[SI]_{2}}{[I]}$, and let $R$ be defined as
\begin{equation}\label{Rlambda}
R = \frac{\tau{w_{1}}{\lambda_{1}}+\tau{w_{2}}{\lambda_{2}}}{\gamma}.
 \end{equation}
Following the method outlined by Keeling \cite{Keeling1999} and Eames \cite{Eames2008}, we
calculate the early quasi-equilibrium values of $\lambda_{1,2}$ as
follows:
\[
\begin{array}{l}
    \dot{\lambda_{1}} = 0 \Leftrightarrow  \dot{[SI]}_{1}[I] = \dot{[I]}[SI]_{1}, \\
    \dot{\lambda_{2}} = 0 \Leftrightarrow  \dot{[SI]}_{2}[I] = \dot{[I]}][SI]_{2}.
 \end{array}
 \]
Upon using the pairwise equations and the closure, consider $[\dot{S}{I}]_{1}{[I]} = [\dot{I}][SI]_{1}$:
\begin{eqnarray}\label{SI1}
   [\dot{S}{I}]_{1}{[I]} &=& ({\tau}w_{1}{[SSI]_{11}}+{\tau}w_{2}{[SSI]_{12}}-{\tau}w_{1}{[ISI]_{11}}-{\tau}w_{2}{[ISI]_{21}}-{\tau}w_{1}{[SI]_{1}}-\gamma{[SI]_{1}})[I]\nonumber\\
   \nonumber \\
    &=& (\tau{w_{1}}{[SI]_{1}}+{\tau}w_{2}{{[SI]_{2}}}-\gamma{[I]})[SI]_{1}.
\end{eqnarray}
Using the classical closure
\[
\begin{array}{l}
 \displaystyle{   \ [ABC]_{12}=\frac{k-1}{k}\frac{[AB]_{1}[BC]_{2}}{[B]} , }\\\\
 \displaystyle{   \ [ABC]_{21}=\frac{k-1}{k}\frac{[AB]_{2}[BC]_{1}}{[B]},}
 \end{array}
 \]
 and making the substitution : $[{S}{I}]_{1}=\lambda_{1}{[I]}$ ,
 $[{S}{I}]_{2}=\lambda_{2}{[I]}$, $[I]\ll1 $, $[S]\approx N$, $[SS]_{1} \approx
 kNp_{1}$, $[SS]_{2} \approx
 kN(1-p_{1})$ together with $\gamma R = \tau{w_{1}}{\lambda_{1}}+\tau{w_{2}}{\lambda_{2}}$, we have
\[
({\tau}w_{1}\lambda_{1}+{\tau}w_{2}\lambda_{2})kp_{1}-({\tau}w_{1}\lambda_{1}+{\tau}w_{2}\lambda_{2})p_{1}-({\tau}w_{1}\lambda_{1}+{\tau}w_{2}\lambda_{2})\lambda_{1}-\tau{w_{1}}\lambda_{1}=0,
\]
which can be solved for $\lambda_1$ to give
\[
\lambda_{1} = \frac{\gamma{(k-1)}{p_{1}}{R}}{\tau{w_{1}}+\gamma{R}}.
\]
Similarly, $\lambda_2$ can be found as
\begin{eqnarray}
\lambda_{2} &=& \frac{\gamma{(k-1)}(1-p_{1}){R}}{\tau{w_{2}}+\gamma{R}}.
\end{eqnarray}
Substituting the expressions for $\lambda_{1,2}$ into the original equation for $R$ yields
\[
R = \frac{A+B+\sqrt{(A+B)^2+4\tau^2w_{1}w_{2}(k-2)}}{2\gamma},
\]
where $A=\tau w_{1}[(k-1)p_{1}-1]$ and $B=\tau w_{2}[(k-1)p_{2}-1]$.
If we define
\[
\displaystyle{R_{1} = \frac{\tau{w_{1}}[(k-1)p_{1}-1]}{\gamma}},\hspace{0.5cm}
\mbox{and}\hspace{0.5cm}\displaystyle{R_{2} = \frac{\tau{w_{2}}[(k-1)p_{2}-1]}{\gamma}},
\]
the expression simplifies to
% \begin{eqnarray}
%       \ R &=& \frac{R_{1}+R_{2}+\sqrt{(R_{1}+R_{2})^2+4R_{1}R_{2} \frac{(k-2)}{((k-1)p_1-1)((k-1)p_2-1))}}}{2},\\
\[
R = \frac{R_{1}+R_{2}+\sqrt{(R_{1}+R_{2})^2+4R_{1}R_{2}Q}}{2},
\]
where $\displaystyle{Q = \frac{(k-2)}{[(k-1)p_1-1][(k-1)p_2-1]}}$.\\\\

\noindent Substituting the modified closure
\[
  \begin{array}{l}
 \displaystyle{  \ [ABC]_{11}=\frac{k_1-1}{k_1}\frac{[AB]_{1}[BC]_{1}}{[B]},}\\\\
 \displaystyle{    \ [ABC]_{12}=\frac{[AB]_{1}[BC]_{2}}{[B]},}\\\\
 \displaystyle{    \ [ABC]_{21}=\frac{[AB]_{2}[BC]_{1}}{[B]},}\\\\
 \displaystyle{    \ [ABC]_{22}=\frac{k_2-1}{k_2}\frac{[AB]_{2}[BC]_{2}}{[B]},}
 \end{array}
 \]
 into (\ref{SI1}) and making further substitution : $[{S}{I}]_{1}=\lambda_{1}{[I]}$, $[{S}{I}]_{2}=\lambda_{2}{[I]}$, $[I]\ll1 $, $[S]\approx N$, $[SS]_{1} \approx
 k_{1} N$, $[SS]_{2} \approx k_{2} N$, we have
\[
 ({\tau}w_{1}\lambda_{1}+{\tau}w_{2}\lambda_{2})k_{1}-({\tau}w_{1}\lambda_{1}+{\tau}w_{2}\lambda_{2})\lambda_{1}-2\tau w_1 \lambda_1 =0\Longrightarrow
 \lambda_{1} = \frac{\gamma{k_{1}}{R}}{2\tau{w_{1}}+\gamma{R}}.
\]
Similarly, the equation $[\dot{S}{I}]_{2}{[I]} = [\dot{I}][SI]_{2}$ yields
\[
\lambda_{2}= \frac{\gamma{k_{2}}{R}}{2\tau{w_{2}}+\gamma{R}}.
\]
Substituting these expressions for $\lambda_{1,2}$ into (\ref{Rlambda}), we have
\[
 \begin{array}{l}
\displaystyle{ R = \frac{\tau{(w_{1}k_{1}+w_{2}k_{2})}-2\tau(w_{1}+w_{2})}{2\gamma}}\\\\
\displaystyle{ \hspace{1cm}+\frac{\sqrt{\left[2\tau(w_{1}+w_{2})-\tau(w_{1}k_{1}+w_{2}k_{2})\right]^2+8\tau^2w_{1}w_{2}(k_1+k_2-2)}}{2\gamma}.}
 \end{array}
\]
If we define
\[
R_{1} = \frac{\tau{w_{1}}(k_{1}-2)}{\gamma},\hspace{0.5cm}R_{2} = \frac{\tau{w_{2}}(k_{2}-2)}{\gamma},
\]
the above expression for $R$ simplifies to
 \begin{equation}
    R = \frac{R_{1}+R_{2}+\sqrt{(R_{1}+R_{2})^2+4R_{1}R_{2}(Q-1)}}{2}
 \end{equation}
 where
 \[
 Q = \frac{k_{1} k_{2}}{(k_{1}-2)(k_{2}-2)}.
 \]

\newpage

%\newpage

\vspace{2cm}

\newpage

\newpage

\newpage

%\newpage
%\begin{figure}
%    \centering
%    \subfigure[$k = 5, \tau = 1$]
%    {
%        \includegraphics[width=7cm]{SIR_k5_10w1.eps}
%        \label{$SIS$}
%    }
%    \subfigure[$k = 10, \tau = 0.5$]
%    {
%        \includegraphics[width=7cm]{SIR_k10_10w1.eps}
%        \label{$SIR$}
%    }
%    \caption{$N = 100, \gamma =
%2, I_{0} = 5 , w_{1} = 10, P(w_{1}) = 0.05, w_{2} = 0.5/0.95$.}
%\end{figure}


\newpage

\begin{thebibliography}{99}
%\addcontentsline{toc}{chapter}{Bibliography}

\bibitem{Al} Almaas, E., Kov\'acs, B., Viczek, T., Oltval, Z.N. \& Barab\'asi, A.-L. (2004). Global organization of metabolic fluxes in the
bacterium {\it Escherichia coli}, {\it Nature} {\bf 427}, 839-843.

\bibitem{AndersonMay} Anderson, R.M. \& May, R.M. (1992). {\it Infectious Diseases of Humans}.  Oxford: Oxford University Press.

\bibitem{Bag} Bagler, G. (2008). Analysis of the airport network of India as a complex weighted network. {\it Physica A} {\bf 387}, 2972-2980.

\bibitem{BallOnKiss} Ball, F. \& Neal, P. (2008). Network epidemic models with two levels of mixing. {\it Math. Biosci.} {\bf 212}, 69-87.

%\bibitem{r15} S. Bansalab, J. Readc, B. Pourbohloulde \& L. A. Meyersfg, The dynamic nature of contact networks in
%infectious disease epidemiology, {\it J. Biological Dynamics} {\bf
%4} 5, 478-489 (2010).

%\bibitem{r2} A. Baronchelli, C. Castellano \& R. Pastor-Satorras 2011. Voter models on weighted networks. Phys. Rev. E.83, 066117.

\bibitem{r1} Barrat, A., Barth\'elemy, M., Pastor-Satorras, R. \& Vespignani, A. (2004). The architecture of complex weighted networks. {\it Proc. Natl. Acad. Sci. USA.} {\bf 101}, 3747-3752.

\bibitem{r7} Barrat, A., Barth\'elemy, M. \& Vespignani, A. (2004). Weighted evolving networks: coupling topology and weight
dynamics. {\it Phys. Rev. Lett.} {\bf 92}, 228701.

\bibitem{r8} Barrat, A., Barth\'elemy, M. \& Vespignani, A. (2004). Modeling the evolution of weighted networks. {\it Phys. Rev. E} {\bf 70},
066149.

\bibitem{r9}  Barrat, A., Barth\'elemy, M. \& Vespignani, A. (2005). The effects of spatial constraints on the evolution of
weighted complex networks. {\it J. Stat. Mech.}, P05003.

\bibitem{BBPV} Barth\'elemy, M., Barrat, A., Pastor-Satorras, R.  \& Vespignani, A. (2005). Characterization and modeling of weighted networks. {\it Physica A} {\bf 346}, 34-43.

\bibitem{BMSKM} Bhattacharya, K., Mukherjee, G., Saram\"aki, J., Kaski, K. \& Manna, S.S. (2008). The international trade network: weighted network analysis and modelling. {\it J. Stat. Mech.}, P02002.

\bibitem{Beu} Beutels, P., Shkedy, Z., Aerts, M. \&  Van Damme, P. (2006). Social mixing patterns for transmission models of
close contact infections: exploring self-evaluation and diary-based data collection through a
web-based interface. {\it Epidemiol. Infect.} {\bf 134}, 1158-1166.

\bibitem{BK05} Blyuss, K.B. \& Kyrychko, Y.N. (2005). On a basic model of a two-disease epidemic. {\it Appl. Math. Comp.} {\bf 160}, 177-187.

\bibitem{BK10} Blyuss, K.B. \& Kyrychko, Y.N. (2010). Stability and bifurcations in an epidemic model with varying immunity period. {\it Bull. Math. Biol.} {\it 72}, 490-505.

\bibitem{Boc} Boccaletti, S., Latora, V., Moreno, Y., Chavez, M. \&  Hwang, D.-U. (2006). Complex networks:
structure and dynamics. {\it Phys. Rep.} {\bf 424}, 175-308.

\bibitem{BDL} Britton, T., Deijfen, M. \& Liljeros, F. (2011). A weighted configuration model and inhomogeneous
epidemics, {\it J. Stat. Phys.} {\bf 145}, 1368-1384.

\bibitem{r13} Britton, T. \& Lindenstrand, D. Inhomogeneous epidemics on weighted networks, {\it Math. Biosci.} published online http://dx.doi.org/10.1016/j.mbs.2012.06.005.

\bibitem{r17} Chu, X., Guan, J., Zhang, Z. \& Zhou, S. (2009). Epidemic spreading in weighted scale-free networks with community
structure. {\it J. Stat. Mech.}, P07043.

\bibitem{Coh97} Cohen, S., Doyle, W.J., Skoner, D.P., Rabin, B.S. \& Gwaltney Jr., J.M. (1997). Social ties and
susceptibility to the common cold, {\it JAMA} {\bf 277}, 1940-1944.

\bibitem{Col} Colizza, V., Barrat, A., Barth\'elemy, M., Valleron, A.-J. \& Vespignani, A. (2007). Modelling the
worldwide spread of pandemic influenza: baseline case and containment interventions. {\it PLoS Med.} {\bf 4}, 95-110.

\bibitem{Col2} Colizza, V., Barrat, A., Barth\'elemy, M. \& Vespignani, A. (2006). The role of the airline transportation network in the prediction and predictability of global epidemics. {\it Proc. Natl. Acad. Sci. USA} {\bf 103}, 2015-2020.

\bibitem{Coop} Cooper, B.S., Pitman, R.J., Edmunds, W.J. \& Gay, N.J. (2006). Delaying the international
spread of pandemic influenza. {\it PLoS Med.} {\bf 3}, e212.

\bibitem{DanonNetwReview} Danon, L., Ford, A.P., House, T., Jewell, C.P., Keeling, M.J., Roberts, G.O., Ross, J.V. \& Vernon, M.C. (2011). Networks
and the epidemiology of infectious disease. {\it Interdisc. Persp. Infect. Diseases} {\bf 2011}, 284909.

%\bibitem{r10} B. Daudert \& Bai-Lian Li, Spreading of infectious disease on complex networks with non-symmetric transmission
%probabilities, published online
%http://trove.nla.gov.au/work/10085638 (2006).

\bibitem{Deij} Deijfen, M. (2011). Epidemics and vaccination on weighted graphs. {\it Math. Biosci.} {\bf 232}, 57-65).

\bibitem{DiHe} Diekmann, O. \& Heesterbeek, J.A.P. (2000). {\it Mathematical epidemiology of infectious diseases: model building, analysis and interpretation}. Chichester: Wiley.

\bibitem{Diekmann1990R0} Diekmann, O., Heesterbeek, J.A.P. \& Metz, J.A.J. (1990). On
the definition and the computation of the basic reproduction ratio $R_0$, in models for infectious diseases in heterogeneous
populations. {\it J. Math. Biol.} {\bf 28}, 365-382.

\bibitem{DM} Dorogovtsev, S.N. \& Mendes, J.F.F. (2003). {\it Evolution of networks: From biological nets to the Internet and WWW}.
Oxford: Oxford University Press.

\bibitem{Eames2008} Eames K.T.D. (2008). Modelling disease spread through random and regular contacts in
clustered populations. {\it Theor. Popul. Biol.} {\bf 73}, 104-111.

\bibitem{HetPairwise} Eames, K.T.D. \& Keeling, M.J. (2002). Modeling dynamic and network heterogeneities in the spread of sexually transmitted
diseases. {\it Proc. Natl. Acad. Sci. USA} {\bf 99}, 13330-13335.

\bibitem{r12} K.T.D Eames, J.M. Read \& W.J. Edmunds, Epidemic prediction and control in weighted networks, {\it Epidemics} {\bf 1}, 70-76 (2009).

\bibitem{Ed97} Edmunds, W.J., O'Callaghan, C.J. \& Nokes, D.J. (1997). Who mixes with whom? A method to
determine the contact patterns of adults that may lead to the spread of airborne infections. {\it Proc. R. Soc. Lond. B}
{\bf 264}, 949-957.

\bibitem{Eubank} Eubank, S., Guclu, H., Kumar, V.S.A., Marathe, M.V., Srinivasan, A., Toroczkai, Z. \& Wang, N. (2004).
Modelling disease outbreak in realistic urban social networks. {\it Nature} {\bf 429}, 180-184.

\bibitem{FRS} Fagiolo, G., Reyes, J. \& Schiavo, S. (2008). On the topological properties of the world trade web: A weighted
network analysis. {\it Physica A} {\bf 387}, 3868-3873.

\bibitem{Fan} Fan, Y., Li, M., Chen, J., Gao, L., Di, Z. \& Wu, J. (2004). Network of econophysicists: a weighted network to investigate the development of econophysics.
{\it Int. J. Mod. Phys. B} {\bf 18}, 17-19, 2505-2512.

\bibitem{r11} Garlaschelli, D. (2009). The weighted random graph model. {\it New J. Phys.}  {\bf 11}, 073005.

\bibitem{GMB} Gilbert, M., Mitchell, A., Bourn, D., Mawdsley, J., Clifton-Hadley, R. \& Wint, W. (2005). Cattle
movements and bovine tuberculosis in Great Britain. {\it Nature} {\bf 435}, 491-496.

\bibitem{Gil} Gillespie, D.T. (1977). Exact stochastic simulation of coupled chemical reactions. {\it J. Phys. Chem.} {\bf 81}, 2340-2361.

\bibitem{GS} Gross, T. \& Sayama, H. (Eds.) (2009). {\it Adaptive networks: theory, models and applications}. New York: Springer.

\bibitem{Vasilis} Hatzopoulos, V., Taylor, M., Simon, P.L. \& Kiss, I.Z. (2011). Multiple sources and routes of information transmission: implications for epidemic
dynamics. {\it Math. Biosci.} {\bf 231}, 197-209.

\bibitem{HouseMotif} House, T., Davies, G., Danon, L. \& Keeling, M.J. (2009). A motif-based approach to network epidemics,
{\it Bull. Math. Biol.} {\bf 71}, 1693-1706.

\bibitem{ThomasUnifying} House, T. \& Keeling, M.J. (2011). Insights from unifying modern approximations to infections on networks. {\it J. Roy. Soc. Interface} {\bf 8},
67-73.

\bibitem{HBG} Hufnagel, L., Brockmann, D. \& Geisel, T. (2004). Forecast and control of epidemics in a globalized world. {\it Proc. Natl. Acad. Sci. USA} {\bf 101}, 15124-15129.

\bibitem{NetwWithSaturation} Joo, J. \& Lebowitz, J.L. (2004). Behavior of susceptible-infected-susceptible epidemics on heterogeneous networks with saturation.
{\it Phys. Rev. E} {\bf 69}, 066105.

\bibitem{Ka} Karsai, M., Juh\'asz, R. \& Ingl\'oi, F. (2006). Nonequilibrium phase transitions and finite-size scaling in weighted scale-free networks. {\it Phys. Rev. E} {\bf 73}, 036116.

\bibitem{Keeling1999} Keeling, M.J. (1999). The effects of local spatial structure on epidemiological invasions. {\it Proc. R. Soc. Lond. B} {\bf 266}, 859-867.

\bibitem{KeRo} Keeling, M.J. \& Rohani, P. (2007). {\it Modeling infectious diseases in humans and animals}. Princeton: Princeton University Press.

\bibitem{NetwReviewKeelingEames} Keeling, M.J. \& Eames, K.T.D. (2005). Networks and epidemic models. {\it J. R. Soc. Interface} {\bf 2}, 295-307.

\bibitem{KissMultipRoutes} Kiss, I.Z., Green, D.M. \& Kao, R.R. (2006). The effect of contact
heterogeneity and multiple routes of transmission on final epidemic size. {\it Math. Biosci.} {\bf 203}, 124-136.

\bibitem{KissInfo} Kiss, I.Z., Cassell, J., Recker, M. \& Simon, P.L. (2010). The impact of information transmission on epidemic outbreaks. {\it Math. Biosci.} {\bf 225}, 1-10.

\bibitem{LiCai} Li, W. \& Cai, X. (2004). Statistical analysis of airport network of China. {\it Phys. Rev. E} {\bf 69}, 046106.

\bibitem{LC} Li, C. \& Chen, G. (2004). A comprehensive weighted evolving network model. {\it Physica A} {\bf 343}, 288-294.

\bibitem{LiWu} Li, M., Wu, J., Wang, D., Zhou, T., Di, Z. \& Fan, Y. (2007). Evolving model of weighted networks inspired by scientific collaboration networks. {\it Physica A} {\bf 375}, 355-364

\bibitem{MoNe} Moore, C. \& Newman, M.E.J. (2000). Epidemics and percolations on small-world networks. {\it Phys. Rev. E} {\bf 61}, 5678-5682.

\bibitem{MPV} Moreno, Y., Pastor-Satorras, R. \& Vespignani, A. (2002). Epidemic outbreaks in complex heterogeneous networks. {\it Eur. Phys. J. B} {\bf 26}, 521-529.

\bibitem{NewEpi} Newman, M.E.J. (2002). Spread of epidemic disease on networks. {\it Phys. Rev. E} {\bf 66}, 016128.

\bibitem{New} Newman, M.E.J. (2004). Analysis of weighted networks. {\it Phys. Rev. E} {\bf 70}, 056131.

\bibitem{NS}  Noda, K., Shinohara, A., Takeda, M., Matsumoto, S., Miyano, S. \& Kuhara, S. (1998). Finding genetic network from experiments by weighted network
model. {\it Gen. Inform.} {\bf 9}, 141-150.

\bibitem{LewiStone} Olinky, R. \& Stone, L. (2004). Unexpected epidemic thresholds in heterogeneous networks: The role of disease transmission.
{\it Phys. Rev. E} {\bf 70}, 030902(R).

\bibitem{OS} Onnela, J.-P., Saram\"aki, J., Hyv\"onen, J., Szab\'o, G., de Menezes, M.A., Kaski, K., Barab\'asi, A.-L. \& Kert\'esz, J. (2007).
Analysis of a large-scale weighted network of one-to-one human communication. {\it New J. Phys.} {\bf 9}, 179.

%\bibitem{PLY} K. Park, Y.-C. Lai \& N. Ye, Characterization of weighted complex networks, {\it Phys. Rev. E} {\bf 70}, 026109 (2004).

\bibitem{PS1} Pastor-Satorras, R. \& Vespignani, A. (2001). Epidemic spreading in scale-free networks. {\it Phys. Rev. Lett.} {\bf 86}, 3200-3202.

\bibitem{PS2} Pastor-Satorras, R. \& Vespignani, A. (2001). Epidemic dynamics and endemic states in complex networks. {\it Phys. Rev. E} {\bf 63}, 066117.

\bibitem{r5} Rand, D.A. (1999). Correlation equations and pair approximations for spatial ecologies. {\it CWI Quarterly} {\bf 12}, 329-368.

\bibitem{Read08} Read, J.M., Eames, K.T.D. \& Edmunds, W.J. (2008). Dynamic social networks and the implications for the spread of infectious
disease. {\it J. R. Soc. Interface} {\bf 5}, 1001-1007.

\bibitem{Riley07} Riley, S. (2007). Large-scale spatial-transmission models of infectious disease. {\it Science} {\bf 316}, 1298-1301.

\bibitem{RFer06} Riley, S. \& Ferguson, N.M. (2006). Smallpox transmission and control: spatial dynamics in
Great Britain. {\it Proc. Natl. Acad. Sci.} {\bf 103}, 12637-12642.

\bibitem{KieranAsymmetric} Sharkey, K.J., Fernandez, C., Morgan, K.L., Peeler, E., Thrush, M., Turnbull, J.F. \& Bowers, R.G. (2006). Pair-level approximations to the spatio-temporal dynamics of epidemics on asymmetric contact networks. {\it J. Math. Biol.} {\bf 53}, 61-85.

%\bibitem{r14} P. Schumm, C. Scoglio, D. Gruenbacher \& T. Easton, Epidemic Spreading on Weighted Contact Networks, {\it
%Bionetics} 10.1109/BIMNICS.4610111 (2007).

%\bibitem{MickMarkovianPairwise} M. Taylor, P.L. Simon, D.M. Green, T. House \& I. Z. Kiss, From Markovian to pairwise epidemic models and the performance of moment closure approximations, {\it %J. Math. Biol.} {\bf 64}, 1021-1042 (2012).

%\bibitem{r6} M. Taylor, P.L. Simon, D.M. Green, T. House \& I.Z. Kiss 2011. From Markovian to pairwise epidemic models and the
%performance of moment closure approximations. J. Math. Biol. 62(4),  479-508.

\bibitem{r3} Wang, N.-N. \& Chen, G.-L. (2011). A virus spread model based on cellular automata in weighted scale-free networks. {\it Am. J. Engrg, Techn. Res.} {\bf 11}, 148-154.

\bibitem{WZ} Wang, S. \& Zhang, C. (2004). Weighted competition scale-free network. {\it Phys. Rev. E} {\bf 70}, 066127.

\bibitem{WW} Wang, W.-X., Wang, B.-H., Hu, B., Yan, G. \& Ou, Q. (2005). General Dynamics of Topology and Traffic onWeighted Technological Networks. {\it Phys. Rev. Lett.} {\bf 94}, 188702.

\bibitem{Webb} Webb, C.R. (2006). Investigating the potential spread of infectious diseases of sheep via agricultural shows in Great Britain. {\it Epidemiol. Infect.} {\bf 134}, 31-40.

\bibitem{Yan} Yan, G., Zhou, T., Wang, J., Fu, Z.-Q. \& Wang, B.-H. (2005). Epidemic spread in weighted scale-free networks. {\it Chinese Phys. Lett.} {\bf 22}, 510.

\bibitem{r16} Yang, Z. \& Zhou, T. (2012). Epidemic spreading in weighted networks: an edge-based mean-field solution. {\it Phys. Rev. E} {\bf 85}, 056106.

\bibitem{Yang} Yang, R., Zhou, T., Xie, Y.-B., Lai, Y.-C. \& Wang, B.-H. (2008). Optimal contact process on complex networks. {\it Phys. Rev. E} {\bf 78}, 066109.

\bibitem{YJBT} Yook, S.H., Jeong, H., Barab\'asi, A.-L. \& Tu, Y. (2001). Weighted evolving networks. {\it Phys. Rev. Lett.} {\bf 86}, 5835-5838.

\bibitem{ZTZH} Zheng, D., Trimper, S., Zheng, B. \& Hui, P.M. (2003). Weighted scale-free networks with stochastic weight assignments. {\it Phys Rev. E} {\bf 67}, 040102R.

\end{thebibliography}
\end{document}